\documentclass[a4paper,12pt]{amsart}

\usepackage{amssymb}

\def\hbar{\bar{h}}

\def\lhd{\vartriangleleft}

\def\iso{\buildrel \sim\over\to}

\def\GS{{\mathfrak{S}}}

\def\Gm{{\mathfrak{m}}}
\def\Gn{{\mathfrak{n}}}
\def\Gp{{\mathfrak{p}}}
\def\Gt{{\mathfrak{t}}}

\def\gl{{\mathfrak{gl}}}
\def\Gsl{{\mathfrak{sl}}}

\def\CA{{\mathcal{A}}}
\def\CB{{\mathcal{B}}}
\def\CC{{\mathcal{C}}}
\def\CD{{\mathcal{D}}}

\def\CH{{\mathcal{H}}}
\def\CI{{\mathcal{I}}}

\def\CO{{\mathcal{O}}}
\def\CP{{\mathcal{P}}}

\def\CS{{\mathcal{S}}}
\def\CT{{\mathcal{T}}}

\def\BA{{\mathbf{A}}}

\def\BC{{\mathbf{C}}}

\def\BG{{\mathbf{G}}}
\def\BH{{\mathbf{H}}}

\def\BL{{\mathbf{L}}}

\def\BP{{\mathbf{P}}}
\def\BQ{{\mathbf{Q}}}
\def\BR{{\mathbf{R}}}

\def\BT{{\mathbf{T}}}

\def\BZ{{\mathbf{Z}}}
\def\Ba{{\mathbf{a}}}

\def\Bc{{\mathbf{c}}}

\def\Bh{{\mathbf{h}}}

\def\Bk{{\mathbf{k}}}

\def\Bq{{\mathbf{q}}}

\def\Bs{{\mathbf{s}}}
\def\Bt{{\mathbf{t}}}

\def\Bx{{\mathbf{x}}}

\def\eps{\varepsilon}

\def\can{{\mathrm{can}}}

\def\coker{\operatorname{coker}\nolimits}

\def\End{\operatorname{End}\nolimits}

\def\Ext{\operatorname{Ext}\nolimits}

\def\GL{\operatorname{GL}\nolimits}

\def\Hilb{\operatorname{Hilb}\nolimits}

\def\Hom{\operatorname{Hom}\nolimits}

\def\Id{\operatorname{Id}\nolimits}
\def\id{\operatorname{id}\nolimits}
\def\im{\operatorname{im}\nolimits}

\def\Ind{\operatorname{Ind}\nolimits}

\def\Irr{\operatorname{Irr}\nolimits}

\def\mMod{\operatorname{\!-mod}\nolimits}

\def\opp{{\operatorname{opp}\nolimits}}

\def\mperf{\operatorname{\!-perf}\nolimits}

\def\Pic{\operatorname{Pic}\nolimits}

\def\mproj{\operatorname{\!-proj}\nolimits}

\def\rank{\operatorname{rank}\nolimits}

\def\Res{\operatorname{Res}\nolimits}

\def\Spec{\operatorname{Spec}\nolimits}

\def\mtilt{\operatorname{\!-tilt}\nolimits}

\def\Tor{\operatorname{Tor}\nolimits}
\def\Tr{\operatorname{Tr}\nolimits}

\def\ie{{\em i.e.}}

\def\tf{{\tilde{f}}}

\def\tDelta{{\tilde{\Delta}}}

\def\tGamma{{\tilde{\Gamma}}}
\def\tnabla{{\tilde{\nabla}}}
\def\tOmega{{\tilde{\Omega}}}

\newtheorem{thm}{Theorem}[section]
\newtheorem{lemma}[thm]{Lemma}
\newtheorem{cor}[thm]{Corollary}
\newtheorem{prop}[thm]{Proposition}
\newtheorem{defi}[thm]{Definition}
\newtheorem{conj}[thm]{Conjecture}

\theoremstyle{definition}
\newtheorem{rem}[thm]{Remark}
\newtheorem{example}[thm]{Example}
\newtheorem{hyp}{Hypothesis}

\input xypic
\xyoption{all}
\def\reg{\operatorname{reg}\nolimits}

\def\KZ{\operatorname{KZ}\nolimits}
\newcommand\CHEVIE{{\sf {CHEVIE}}}
\newcommand\GAP{{\sf {GAP}}}

\begin{document}
\author{Rapha\"el Rouquier}
\address{Rapha\"el Rouquier\\
Mathematical Institute\\
24-29 St Giles'\\
Oxford, OX1 3LB, UK}
\email{rouquier@maths.ox.ac.uk}
\title{$q$-Schur algebras and complex reflection groups}
\dedicatory{To Victor Ginzburg, on his fiftieth birthday}

\begin{abstract}
We show that the category $\CO$ for a rational Cherednik algebra of type $A$
is equivalent
to modules over a $q$-Schur algebra (parameter $\not\in\frac{1}{2}+\BZ$),
providing thus character formulas for simple modules. We give some
generalization to $B_n(d)$. We prove an ``abstract'' translation
principle. These results follow from the unicity of certain highest
weight categories
covering Hecke algebras. We also provide a semi-simplicity criterion for Hecke
algebras of complex reflection groups and show the isomorphism type of
Hecke algebras is invariant under field automorphisms acting on parameters.
\end{abstract}

\maketitle
\section{Introduction}
This paper (and its sequel) develops a new aspect of the representation
theory of Hecke algebras of complex reflection groups, namely the study of
quasi-hereditary covers, analogous to $q$-Schur algebras
in the symmetric groups case. An important point is the existence of
a family of such covers: it depends on the choice of ``logarithms'' of
the parameters.

\smallskip
The theory we develop is particularly interesting when the ring of
coefficients is not specialized: it blends features of representation
theory over $\BC$ at roots of unity and features away from roots of unity,
where Lusztig's families of characters show up (in that respect, it
is a continuation of \cite{Rou1}, where combinatorial objects are given
homological definitions, which led to generalizations from real to
complex reflection groups).

\smallskip
The main idea of this first paper is the unicity of certain types of
quasi-hereditary covers. This applies in particular to the category $\CO$
of rational Cherednik algebras~: we show that, in type $A$, when the parameter is not
in $\frac{1}{2}+\BZ$, the category $\CO$ is equivalent to the module category
of a $q$-Schur algebra, solving a conjecture of \cite{GGOR}. As a
consequence, we obtain character formulas for simple objects of $\CO$ in
this case. We also obtain a general translation principle for category $\CO$
of a Cherednik algebra.

\medskip
In \S \ref{sectioncplx}, we introduce a function ``c'' on the set
of irreducible characters of $W$, with values
linear functions of the logarithms of the parameters and we construct
an order on the set of irreducible characters of $W$. This is suggested
by \cite[Lemma 2.5]{DunOp} (``roots of unity'' case) as well as by \cite{Lu1}
(``away from roots of unity'').

In \S \ref{sectionqhcovers}, we develop a general theory of (split)
highest weight categories over a commutative ring. This is a categorical
version of Cline-Parshall-Scott's integral quasi-hereditary algebras.
We study covers of finite dimensional algebras by highest weight categories
and consider different levels of ``faithfulness''. The simplest
situation is that of a ``double centralizer Theorem''. The key results are
Proposition \ref{equivDelta} (deformation principle) and
Theorem \ref{unicityfromorder} (unicity).

\S \ref{sectionCherednik} shows that category $\CO$ for Cherednik
algebras gives a cover of Hecke algebras of complex reflection groups, and that
it has the faithfulness property when the rank $1$ parabolic Hecke
subalgebras are
semi-simple. This provides a translation principle for category $\CO$.
We also give a simple criterion for semi-simplicity of Hecke algebras
in characteristic $0$, generalizing the usual property for Coxeter groups
and the ``one-parameter case'', that the algebra is semi-simple if that
parameter is not a root of unity (Theorem \ref{ss}). We prove that
rescaling the parameters by a positive integer (without affecting
``denominators'') doesn't change category $\CO$, up to equivalence, 
and that the Hecke algebra is unchanged, up to
isomorphism of $\BC$-algebras, by field automorphisms acting on parameters.
In the last part, we describe blocks with ``defect $1$''.

Finally, in \S \ref{sectionBnd}, we consider the case $W=B_n(d)$.
We show that, for a suitable choice of ``logarithms of parameters'', the
category $\CO$ is equivalent to modules over Dipper-James-Mathas
$q$-Schur algebra (Theorem \ref{compBnd}). Otherwise, we obtain new $q$-Schur
algebras. Putting our work together with Yvonne's \cite{Yv} suggests that the
decomposition matrices should be given by Uglov's canonical
bases of the level $d$ Fock spaces.

\medskip
Relations between Kazhdan-Lusztig theory and modular representations
at roots of unity have been investigated by various authors
\cite{Ge2,GeRou,Ge3,Jac1,Jac2}, and \cite{DuPaSc1,DuPaSc2,DuPaSc3}, whose
``integral''
approach influenced our \S \ref{sectionqhcovers}. We hope our approach
provides some new insight.

The second part will deal with integral aspects, bad primes and dualities
and will address the relations between the representation theory
``at $t=0$'' of the rational Cherednik algebra and Lusztig's
asymptotic Hecke algebra.
We will discuss more thoroughly the case of finite Coxeter
groups and present a number of conjectures.

\bigskip
I thank Susumu Ariki, Steve Donkin,
Karin Erdmann, Pavel Etingof, Victor Ginzburg and Bernard Leclerc for useful
discussions.

\section{Notations}
Let $k$ be a commutative ring and $A$ a $k$-algebra. We denote by
$A\mMod$ the category of finitely generated $A$-modules and by
$A\mproj$ the category of finitely generated projective $A$-modules.
We write $\otimes$ for $\otimes_k$.
Let $M$ be a $k$-module.
We put $M^*=\Hom_k(M,k)$ and given $n$ a non-negative integer, we write
$M^{\oplus n}$ for $M^n$, the direct sum of $n$ copies of $M$, when there is
a risk of confusion.

Let $k'$ be a commutative $k$-algebra.
We put
$k'M=k'\otimes M$. We put $k'(A\mMod)=(k'A)\mMod$. We denote
by $\Irr(A)$ the set of isomorphism classes of simple $A$-modules.
If $\Gm$ is a maximal ideal of $k$, then we put $M(\Gm)=(k/\Gm)M$, etc...
Given $B$ another $k$-algebra, we write $(A\mMod)\otimes (B\mMod)$ for
$(A\otimes B)\mMod$.

\smallskip
Let $\CA$ be an abelian category.
We denote by $D^b(\CA)$ the derived category of bounded complexes of objects
of $\CA$.
We denote by $\CA\mproj$ the full subcategory of $\CA$ of projective objects.
Given $I$ a set of objects of $\CA$,
we denote by $\CA^I$ the full exact subcategory of $\CA$ of
$I$-filtered objects, \ie, objects
that have a finite filtration whose successive quotients are isomorphic
to objects of $I$.

Given $G$ a finite group, we denote by $\Irr(G)$ the set of irreducible
(complex) characters of $G$. 

Let $\Lambda$ be a set. Given $\le_1$ and $\le_2$ two orders on
$\Lambda$, we say that $\le_1$ refines $\le_2$ if 
$\lambda\le_2 \lambda'$ implies $\lambda\le_1 \lambda'$.
Fix an order on $\Lambda$. A subset $I$ of $\Lambda$ is an ideal (resp.
a coideal)
if $\lambda'\le \lambda$ (resp. $\lambda\le \lambda'$) and $\lambda\in I$
imply $\lambda'\in I$. Given $\lambda\in\Lambda$, we put
$\Lambda_{<\lambda}=\{\lambda'\in\Lambda | \lambda'<\lambda\}$, etc...

\section{Parameters for Hecke algebras}
\label{sectioncplx}

\subsection{Definitions}
\subsubsection{Hecke algebra}
\label{secHecke}
Let $W$ be a finite reflection group on a complex vector space $V$.
Let $\CA$ be the set of reflecting hyperplanes of $W$ and for $H\in\CA$,
let $W_H$ be the pointwise stabilizer of $H$ in $W$, let $e_H=|W_H|$, and
let $o_H$ be the cardinality of $W(H)$ (=orbit of $H$ under $W$).

Let $U=\coprod_{H\in\CA/W}\Irr(W_H)$.
We have a bijection $\BZ/e_H\iso \Irr(W_H),\ j\mapsto \det^j_{|W_H}$,
and we denote by $(H,j)$ the corresponding element of $U$.
Let $\BG_m$ be the multiplicative group over $\BZ$.
Let $\BT=(\BG_m)^U$ and $\Bk=\BZ[\BT]=\BZ[\{\Bx_u^{\pm 1}\}_{u\in U}]$.

\smallskip
Let $V_{\reg}=V-\bigcup_{H\in\CA}H$, let $x_0\in V_{\reg}$, and let
$B_W=\pi_1(V_{\reg}/W,x_0)$ be the braid group of $W$.
Let $\BH$ be the Hecke algebra of $W$ over $\Bk$ \cite[\S 4.C]{BrMaRou},
quotient of $\Bk[B_W]$ by the relations
$$\prod_{0\le j< e_H}(\sigma_H-\Bx_{H,j})=0.$$
There is one such relation for each $H\in\CA$. Here, 
$\sigma_H$ is an $s_H$-generator of the monodromy around $H$, where
$s_H$ is the reflection around $H$ with determinant $e^{2i\pi/e_H}$.

\smallskip
In the rest of the paper, we make the following assumption, which is 
known to hold for all but finitely many irreducible complex reflection
groups, for which it is conjectured to be true \cite[\S 4.C]{BrMaRou}
(cf \cite[\S 6]{EtRa} for a proof of a weak version of the conjecture, when
$\dim V=2$).
\begin{hyp}
\label{hyp}
The algebra $\BH$ is free over $\Bk$, of rank $|W|$.
\end{hyp}

\subsubsection{Specialization}
\label{spec}
Let $k$ be a commutative ring.
A {\em parameter} for $W$ is an element $x_\cdot=\{x_u\}$ of $\BT(k)$. This is
the same data as a morphism of groups $X(\BT)\to k^\times$ or
a morphism of rings $\Bk\to k,\ \Bx_u\mapsto x_u$.

\smallskip
Let $m=\mathrm{lcm}(\{e_H\}_{H\in\CA})$ and $\Phi_m(t)\in\BZ[t]$ be the
$m$-th cyclotomic polynomial. Let $k_m=\BZ[t]/(\Phi_m(t))$.
We will identify $k_m$ with its image through the embedding
$k_m\to\BC,\ t\mapsto e^{2i\pi/m}$.

The canonical morphism $\Bk[B_W]\to \Bk[W],\ \sigma_H\mapsto s_H$,
 induces an isomorphism
$k_m\otimes_{\Bk}\BH\iso k_m[W]$ where
the specialization $\Bk\to k_m$ is given by
$\rho=\{t^{jm/e_H}\}_{(H,j)\in U}\in\BT(k_m)$.
\smallskip
It is convenient to shift the elements of $\BT$ by $\rho$.
We put $\Bq_{H,j}=\Bx_{H,j} t^{-jm/e_H}\in k_m[\BT]$. Given a specialization
$k_m[\BT]\to k$, we denote by $q_\cdot$ the image of $\Bq_\cdot$. We put
$\BH(q_\cdot)=k\BH$.

\smallskip
The algebra $k_m\BH$ is a deformation of $k_m[W]$. It follows that
$\BC(\BT)\BH$ is semi-simple.
Let $K$ be a field extension of
$\BC(\BT)$ such that $K\BH$ is split semi-simple. Let $S$ be a local
$\BC[\BT]$-subalgebra of $K$, integrally closed in $K$, and whose
maximal ideal contains $\{\Bq_u-1\}_{u\in U}$.
Then we have a
canonical isomorphism $\Irr(W)\iso\Irr(K\BH)$ (``Tits deformation
Theorem''). More generally, let $k$ be a field such that
$k\BH$ is split semi-simple, together with an integrally closed local
$k_m[\frac{1}{|W|}][\BT]$-subalgebra $S$ of $k$ whose maximal ideal
contains $\{\Bq_u-1\}_{u\in U}$.
Then we have a canonical isomorphism
$\Irr(W)\iso\Irr(k\BH),\chi\mapsto \chi_k$.

\smallskip
By \cite[Corollary 4.8]{Mal}, if the representation $V$ of $W$
is defined over a subfield $K_0$ of $\BC$
and the group of roots of unity in $K_0$ is
finite of order $l$, then $K_0(\{\Bq_u^{1/l}\}_{u\in U})$ is a splitting field
for $\BH$. We choose $S=K_0[\{\Bq_u^{1/l}\}]_{(\Bq_u^{1/l}-1)_u}$
to define the
bijection $\Irr(W)\iso\Irr(K_0(\{\Bq_u^{1/l}\}_{u\in U})\BH)$.

\begin{rem}
One could also work with the smaller coefficient ring 
$\BZ[\{\Ba_u\}][\{\Ba_{H,0}^{\pm 1}\}]$ instead of $\Bk$ and define $\BH$ with
the
relations $\sigma_H^{e_H}+\Ba_{H,e_{H-1}}\sigma_H^{e_H-1}+\cdots+\Ba_{H,0}=0$.
\end{rem}

\subsection{Logarithms of the parameters}
\subsubsection{Function $c$}
Let $\Gt$ be the Lie algebra of $\BT$ over $k_m$.
Let $\{\Bh_u\}_{u\in U}$ be the basis of $X(\Gt)$ giving
the isomorphism $\sum_u \Bh_u:\Gt\iso k_m^U$
corresponding to the isomorphism
$\sum_u \Bq_u:\BT(k_m)\iso \BG_m(k_m)^U$. We denote by
$\Gt_\BZ=\bigoplus_u \Bh_u^{-1}(\BZ)$ the corresponding $\BZ$-Lie
subalgebra of $\Gt$.

\smallskip
Let $\chi\in\Irr(W)$. We put
$$n_{H,j}(\chi)=
\frac{o_He_H}{\chi(1)}\langle\chi_{|W_H},{\det}_{|W_H}^{j} \rangle.$$
This is the scalar by which
$\sum_{H'\in W(H), w\in W_{H'}}\det(w)^{-j} w$ acts on an irreducible
representation of $W$ with character $\chi$.
In particular, this is a non-negative integer.

We define a map $\Bc:\Irr(W)\to X(\Gt)$ by
$$\chi\mapsto \Bc_\chi= \sum_{(H,j)\in U}n_{H,j}(\chi)\Bh_{H,j}.$$

We also put
$$\Bc'_\chi=
\sum_{(H,j)\in U}n_{H,j}(\chi)(\Bh_{H,j}-\Bh_{H,0})=
\Bc_\chi-\sum_{H\in\CA/W} o_He_H\Bh_{H,0}.$$
So, $\Bc'_\chi=0$ if and only if $\chi$ is the trivial character.

\subsubsection{Lift}
\label{choice}
Let $k$ be a commutative ring and $q_\cdot\in (k^\times)^U$. Let
$\Gamma$ be the subgroup of $k^\times$ generated by $\{q_u\}_{u\in U}$.
We denote by $\Gamma_{tor}$ its subgroup of elements of finite order.

Let $\tGamma$ be a free abelian group together with a surjective
morphism $\exp:\tGamma\to \Gamma$ and an isomorphism
$\BZ\iso\ker(\exp)$:
$$0\to \BZ\to \tGamma\xrightarrow{\exp} \Gamma\to 0.$$

Let us fix an order on $\tGamma$ with the following properties:
\begin{itemize}
\item it extends the natural order on $\BZ$
\item it is compatible with the group law
\item if $x\notin\BZ$ and $x>0$, then $x+n>0$ for all $n\in \BZ$, \ie,
the order on $\tGamma$ induces an order on $\Gamma$.
\end{itemize}

We define the {\em coarsest} order to be the one given by
$x>0$ if and only if $x\in\BZ_{>0}$.

\smallskip

Let $h_\cdot\in \tGamma\otimes_\BZ\Gt_\BZ$ with $q_\cdot=\exp(h_\cdot)$~:
this is the data of
$\{h_u\}\in \tGamma^U$ with $q_u=\exp(h_u)$. To $h_\cdot$ corresponds a morphism
$X(\Gt)\to\tGamma$.
We denote by $c:\Irr(W)\to\tGamma$ the map deduced from $\Bc$.

\smallskip
Let $\pi\in B_W$ be the element given by the loop
$t\in [0,1]\mapsto e^{2i\pi t}$. This is a central element of $B_W$
and we denote by $T_\pi$ its image in $\BH$.

We have $\chi_{k}(T_\pi)=\exp(c_\chi)$ \cite[Proposition 4.16]{BrMi}.
Cf also \cite[\S 1]{BrMaMi} for a more detailed discussion.

\begin{rem}
Note that given $\Gamma$, there exists
always $\tilde{\Gamma}$ as above, when $k$ is a domain:
take $\Gamma=\Gamma_{tor}\times L$ with $L$ free and $g$ a generator of
$\Gamma_{tor}$. Let $\tilde{\Gamma}=\langle \tilde{g}\rangle\times L$.
Define
$\exp$ by $\tilde{g}\mapsto g$ and as the identity on $L$. The coarsest
order on $\tilde{\Gamma}$ satisfies the conditions above.
\end{rem}

\begin{example}
Assume the $q_u$'s are roots of unity and $k$ is a domain. Then,
$\Gamma$ is a finite cyclic group and $\tGamma$ is free of rank $1$.
The order on $\tGamma$ is the coarsest order.
\end{example}

\subsubsection{Order on $\Irr(W)$}
\label{order}
We define now an order on $\Irr(W)$.
Let $\chi,\chi'\in\Irr(W)$.
We put $\chi>\chi'$ if $c_\chi<c_{\chi'}$ (equivalently, $c'_\chi<c'_{\chi'}$).

\subsection{Change of parameters}
\subsubsection{Twist}
Let $W^\wedge=\Hom(W,\BC^\times)$ be the group of one-dimensional characters
of $W$. We have an isomorphism given by restriction
$W^\wedge\iso\left(\prod_{H\in\CA}\Irr(W_H)\right)/W$.
The group $W^\wedge$ acts by multiplication on $U$, and this gives an action on
$\Bk_m=k_m\otimes_\BZ \Bk$.
Let $\xi\in W^\wedge$. The action of $\xi$ on $\Bk_m$ is given
by $\Bq_{H,j}\mapsto \Bq_{H,j+r_H}$, where $\xi_{|W_H}=\det^{r_H}_{|W_H}$.
It extends to an automorphism of $k_m$-algebras
$\Bk_m[B_W]\iso \Bk_m[B_W],\ \sigma_H\mapsto \xi(s_H)^{-1}\sigma_H$.
It induces
automorphisms of $k_m$-algebras $k_m\BH\iso k_m\BH$ and
$k_m[W]\iso k_m[W],\ w\mapsto \xi(w)^{-1}w$ for $w\in W$.

There is a similar action of $\xi$ on $X(\Gt)$ given by
$\Bh_{H,j}\mapsto \Bh_{H,j+r_H}$.
We denote by $\theta_\xi$ these automorphisms induced by $\xi$.
We have $\theta_\xi(\Bc_\chi)=\Bc_{\chi\otimes\xi}$.

\subsubsection{Permutation of the parameters}
Consider $G=\prod_{H\in\CA/W}\GS(\Irr(W_H))\subset\GS(U)$. It acts
on $\BT$, hence on $\Bk$. Let $g\in G$. We denote by $\Bk_g$ the
ring $\Bk$ viewed as a $\Bk$-module by letting $a\in\Bk$ act by
multiplication by $g(a)$. There is an isomorphism of
$\Bk$-algebras $\BH\iso\Bk_g \BH,\ \sigma_H\mapsto \sigma_H, a\mapsto g(a)$ for
$a\in \Bk$.
The action on $\Bt$ is given by $h_{H,j}\mapsto h_{H,g(j)}+
\frac{g(j)-j}{e_H}$ (we view $g$ as an automorphism of $\{0,\ldots,e_H-1\}$).

Let $K=\BC(\{\Bq_u^{1/l}\}_u)$.
We extend the action of $G$ to an action by $\BC$-algebra automorphisms
on $K$~:
the element $g$ sends $\Bq_{H,j}^{1/l}$ to
$\Bq_{H,g(j)}^{1/l}e^{2i\pi(g(j)-j)/(le_H)}$. We deduce an action (by
ring automorphisms) on $K\BH$ fixing the image of $B_W$.
The action of $G$ on $\Irr(K\BH)$ induces an action on $\Irr(W)$.

\subsubsection{Normalization}
Consider a map $f:\CA/W\to k^\times$. Let $q'_\cdot$ be given by
$q'_{H,j}=f(H) q_{H,j}$. Let $k'$ be $k$ as a ring, but viewed as a
$\Bk$-algebra through $q'_\cdot$.
Then, we have an isomorphism of $k$-algebras
$k\BH\iso k'\BH,\ T_H\mapsto f(H)^{-1} T_H$ (here, $T_H$ is the image
of $\sigma_H$).

So, up to isomorphism, $k\BH$ depends only on $q_\cdot$ modulo the ``diagonal''
subgroup $(\BG_m)^{\CA/W}$ of $\BT$.
In particular every Hecke algebra over $k$ is isomorphic to one
where $q_{H,0}=1$ for all $H\in\CA$.

\smallskip
Similarly, consider a map $\tilde{f}:\CA/W\to\tGamma$. Put
$h'_{H,j}=\tf(H)+h_{H,j}$. Then, 
${\Bc'_\chi}_{|\Bh_\cdot=h'_\cdot}={\Bc'_\chi}_{|\Bh_\cdot=h_\cdot}$~:
we can reduce to the study of the order on $\Irr(W)$ to the
case where $h_{H,0}=0$ for all $H$.

\begin{rem}
Assume there is $\kappa\in \tGamma$ with
$h_{H,j}=0$ for all $H$ and $j\not=0$ and
$h_{H,0}=\kappa$ (``spetsial case'').
Then, $c'_\chi=-\frac{N(\chi)+N(\chi^*)}{\chi(1)}\kappa$, where $N(\chi)$ is the
derivative at $1$ of the fake degree of $\chi$ (cf \cite[\S 4.B]{BrMi}).

Assume furthermore that $W$ is a Coxeter group. Let
$a_\chi$ (resp. $A_\chi$) be the valuation (resp. the degree)
of the generic degree of $\chi$.
Then, $\frac{N(\chi)+N(\chi^*)}{\chi(1)}=a_\chi+A_\chi$ \cite[4.21]{BrMi},
hence $c'_\chi=-(a_\chi+A_\chi)\kappa$.
\end{rem}

\subsection{Semi-simplicity}
Let us close this part with a semi-simplicity criterion for Hecke
algebras of complex reflection groups over a field of characteristic $0$.
It generalizes the classical idea for Coxeter groups, that, in the equal
parameters case
($(q_{H,0},q_{H,1})=(q,-1)$), the Hecke algebra is semi-simple
if $q$ is not a non-trivial root of unity.

\begin{thm}
\label{ss}
Let $\Bk\to k$ be a specialization with $k$ a characteristic $0$ field.
Assume that $\Gamma_{tor}=1$.
Then, $k\BH$ is semi-simple.
\end{thm}

The proof uses rational Cherednik algebras and will be given in
\S \ref{KZ}.

\medskip
In general, using the action of $T_\pi$ (cf \S \ref{choice}),
we have the following weaker statement:

\begin{prop}
Let $R$ be a local commutative noetherian $\Bk$-algebra with
field of fractions $K$ and residue field $k$. Assume
$K\BH$ is split semi-simple.

If $\chi_K$ and $\chi'_K$ are in the same block of $R\BH$, then
$c_\chi-c_{\chi'}\in\BZ$.
\end{prop}

\section{Quasi-hereditary covers}
\label{sectionqhcovers}
\subsection{Integral highest weight categories}
In this part, we define and study highest weight categories over a
commutative noetherian ring (extending the classical notion for a field).
This matches the definition of split quasi-hereditary algebras \cite{CPS2}.
In the case of a local base ring, a different but equivalent
approach is given in
\cite[\S 2]{DuSc1} (cf also {\em loc. cit.} for comments on the general case).

\subsubsection{Reminders}
Let $k$ be a commutative noetherian ring.
Let $A$ be a finite projective $k$-algebra (\ie, a $k$-algebra,
finitely generated and projective as a $k$-module). Let $\CC=A\mMod$.

\smallskip
Let us recall some basics facts about projectivity.

\smallskip
\noindent
Let $M$ be a finitely generated $k$-module.
The following assertions are equivalent:
\begin{itemize}
\item
$M$ is a projective $k$-module.
\item $k_\Gm M$ is a projective
$k_\Gm$-module for every maximal ideal $\Gm$ of $k$.
\item $\Tor_1^k(k/\Gm,M)=0$ for every maximal ideal $\Gm$ of $k$.
\end{itemize}
Let $M$ be a finitely generated $A$-module. The following assertions are
equivalent:
\begin{itemize}
\item
$M$ is a projective $A$-module.
\item
$k_\Gm M$ is a projective $k_\Gm A$-module for every maximal ideal $\Gm$ of
$k$.
\item
$M$ is a projective $k$-module and
$M(\Gm)$ is a projective
$A(\Gm)$-module for every maximal ideal $\Gm$ of $k$.
\item
$M$ is a projective $k$-module and
$\Ext^1_A(M,N)=0$ for all $N\in \CC\cap k\mproj$.
\end{itemize}

We say that a finitely generated $A$-module $M$ is
{\em relatively $k$-injective} if it is a projective $k$-module and
$\Ext^1_\CC(N,M)=0$ for all $N\in \CC\cap k\mproj$. So, $M$ is
relatively $k$-injective if and only if $M$ is a projective $k$-module and
$M^*$ is a projective right $A$-module.

\subsubsection{Heredity ideals and associated modules}
\label{defqhc}

\begin{defi}
An ideal $J$ of $A$ is an {\em indecomposable split heredity ideal} 
\cite[Definition 3.1]{CPS2} if the following conditions
hold
\begin{itemize}
\item[(i)] $A/J$ is projective as a $k$-module
\item[(ii)] $J$ is projective as a left $A$-module
\item[(iii)] $J^2=J$
\item[(iv)] $\End_A(J)$ is Morita equivalent to $k$.
\end{itemize}
\end{defi}

\begin{rem}
Note that a split heredity ideal, as defined in \cite[Definition 3.1]{CPS2},
is a direct sum of indecomposable split heredity ideals, corresponding
to the decomposition of the endomorphism ring
into a product of indecomposable algebras. Note further that $J$ is
a split heredity ideal for $A$ if and only if it is a split heredity ideal
for the opposite algebra $A^\opp$ \cite[Corollary 3.4]{CPS2}.
\end{rem}

\smallskip
Given $L$ an $A$-module, we denote by
$$\tau_L:L\otimes_{\End_A(L)} \Hom_A(L,A)\to A, l\otimes f\mapsto f(l)$$
the canonical morphism of $(A,A)$-bimodules.

Given $P$ an $A$-module, 
we define similarly $\tau'_{L,P}:L\otimes\Hom_A(L,P)\to P$.

\begin{lemma}
\label{idempotent}
Let $L$ be a projective $A$-module. 
Then, $J=\im\tau_L$ is an ideal of $A$ and $J^2=J$.
\end{lemma}

\begin{proof}
Since $\tau_L$ is a morphism of $(A,A)$-bimodules, it is clear that
$J$ is an ideal of $A$. Let $E=\End_A(L)$ and $L^\vee=\Hom_A(L,A)$.
We have a commutative diagram
$$\xymatrix{
L\otimes_E L^\vee\otimes_A L\otimes_E L^\vee\ar[dr]_{\tau_L\otimes\tau_L}
\ar[rr]_\sim^{l\otimes f\otimes l'\otimes f'\mapsto l\otimes f(-)l'\otimes f'}
 &&
L\otimes_E \End_A(L)\otimes_E L^\vee\ar[dl]^{\ \ \ \ \ \
l\otimes\phi\otimes f'\mapsto f'(\phi(l))} \\
& A\otimes_A A=A
}$$
where the horizontal arrow is an isomorphism since $L$ is projective.
The image in $A$ of $L\otimes_E \id_L\otimes_E L^\vee$ is equal to $J$ and
the diagram shows it is contained in $J^2$.
\end{proof}

\begin{lemma}
\label{orthogonal}
Let $J$ be an ideal of $A$ which is projective as a left $A$-module and such
that $J^2=J$.

Let $M$ be an $A$-module. Then, 
$\Hom_A(J,M)=0$ if and only if $JM=0$.
\end{lemma}

\begin{proof}
Consider $m\in M$ with $Jm\not=0$. The morphism of $A$-modules
$J\to M,\ j\mapsto jm$ is not zero. This shows the first implication.
The image of a morphism of $A$-modules $J\to M$ is
contained in $JM$, since $J^2=J$. This proves the Lemma.
\end{proof}

\begin{lemma}
\label{hermod}
Let $L$ be a projective object of $\CC$ which is a faithful $k$-module.
The following assertions are equivalent
\begin{itemize}
\item[(i)]
$\tau'_{L,P}:L\otimes \Hom_\CC(L,P)\to P$
is a split injection of $k$-modules for
all projective objects $P$ of $\CC$.
\item[(ii)]
$\tau'_{L,A}:L\otimes \Hom_\CC(L,A)\to A$ is a split injection of $k$-modules.
\item[(iii)]
$k\iso\End_\CC(L)$ and
given $P$ a projective object of $\CC$, then 
there is a subobject $P_0$ of $P$ such that
\begin{itemize}
\item $P/P_0$ is a projective $k$-module
\item $\Hom_\CC(L,P/P_0)=0$ and
\item
$P_0\simeq L\otimes U$ for some $U\in k\mproj$.
\end{itemize}
\end{itemize}
\end{lemma}

\begin{proof}
The equivalence between (i) and (ii) is clear.

Assume (i). Then, $\tau'_{L,L}:L\otimes\End_A(L)\to L$ is injective. Since
it is clearly surjective, it is an isomorphism. Since $L$ is a progenerator
for $k$, we obtain $k\iso\End_A(L)$. Let $P$ be a projective object of
$\CC$. Let $P_0=\im\tau'_{L,P}$, a direct summand of $P$ as a $k$-module.
The map
$$\Hom_A(L,\tau'_{L,P}):\Hom_A(L,L\otimes \Hom_A(L,P))\to\Hom_A(L,P)$$
is clearly surjective, hence $\Hom_A(L,P/P_0)=0$. This proves (iii).

Assume (iii). Let $P$ be a projective object of $\CC$.
The canonical map $\Hom_A(L,P_0)\to\Hom_A(L,P)$ is an
isomorphism.
We have canonical isomorphisms
$\Hom_A(L,L\otimes U)\iso \End_A(L)\otimes U\iso U$.
So, $\tau'_{L,P}$ is injective with image $P_0$ and (i) holds.
\end{proof}

\begin{rem}
Note that if $k$ has no non-trivial idempotent, then every non-zero
projective $k$-module is faithful.
\end{rem}

Let $M(\CC)$ be the set of isomorphism classes of projective
objects $L$ of $\CC$ satisfying the equivalent assertions of Lemma
\ref{hermod}.

Let $\Pic(k)$ be the group of isomorphism classes of invertible
$k$-modules.
Given $F\in\Pic(k)$ and
$L\in M(\CC)$, then $L\otimes F\in M(\CC)$. This gives an action of
$\Pic(k)$ on $M(\CC)$.

\begin{prop}
\label{moduleideal}
There is a bijection from $M(\CC)/\Pic(k)$ to
the set of indecomposable split heredity ideals of $A$ given by
$L\mapsto \im(\tau_L)$.

Furthermore, the canonical functor $(A/\im\tau_L)\mMod\to A\mMod$
induces an equivalence between $(A/\im\tau_L)\mMod$ and the full subcategory
of $\CC$ of objects $M$ such that $\Hom_\CC(L,M)=0$.
\end{prop}

\begin{proof}
We will prove a more precise statement.
We will construct inverse maps $\alpha,\beta$ between $M(\CC)$ and the set of 
isomorphism classes of pairs
$(J,P)$, where $J$ is an indecomposable split heredity ideal of $A$ and
$P$ is a progenerator for $\End_A(J)$ such that $k\iso\End_{\End_A(J)}(P)$.
Here, we say that two pairs $(J,P)$ and $(J',P')$ are isomorphic if
$J'=J$ and $P'\simeq P$.

\smallskip
Let $L\in M(\CC)$, let $J=\im\tau_L$ and let $B=\End_A(J)$. By assumption,
$A/J$ is a projective $k$-module. Also, $J^2=J$ by
Lemma \ref{idempotent}. 
Note that $\Hom_A(L,A)$ is a faithful projective $k$-module.
Since $L\otimes \Hom_A(L,A)\simeq J$, it follows that $J$ is a projective
$A$-module. Also, $\End_k(\Hom_A(L,A))\iso\End_A(J)$ because
$k\iso\End_A(L)$. This gives $\Hom_A(L,A)$ a structure of right $B$-module.
Let $P=\Hom_k(\Hom_A(L,A),k)$. This is a progenerator for $B$ and
$k\iso\End_B(P)$. We have obtained a pair
$(J,P)=\alpha(L)$ as required.

\smallskip
Consider now a pair $(J,P)$. Let $B=\End_A(J)$.
Let $L=J\otimes_B P$, a projective $A$-module. We have
$k\iso\End_B(P)\iso\End_A(L)$.
Let $i:J\to A$ be the inclusion map. There is $p\in P$ and
$f\in\Hom_B(P,\Hom_A(J,A))$ such that $f(p)=i$. Let
$g\in\Hom_A(J\otimes_B P,A)$ be the adjoint map.
Given $j\in J$, we have
$\tau_{J\otimes_B P}(j\otimes_B p\otimes g)=j$. So,
$J\subset\im\tau_L$.
Finally,
we have $\Hom_A(L,A/J)\iso\Hom_B(P,\Hom_A(J,A/J))$.
By Lemma \ref{orthogonal},
we have $\Hom_A(J,A/J)=0$, hence $\Hom_A(L,A/J)=0$. So,
$\im(\tau_L)\subset J$, hence $\im(\tau_L)=J$.
We have an isomorphism of right $B$-modules
$\Hom_B(P,B)\simeq \Hom_k(P,k)$ by Morita theory.
We have $\End_A(J)\iso\Hom_A(J,A)$ since $\Hom_A(J,A/J)=0$.
So, we have isomorphism of right $B$-modules
$$\Hom_A(L,A)\simeq \Hom_B(P,\Hom_A(J,A))\simeq \Hom_B(P,B)\simeq
\Hom_k(P,k).$$
We deduce 
$$J\simeq J\otimes_B P\otimes \Hom_k(P,k)\iso L\otimes \Hom_A(L,A).$$
Now, $\tau_L:L\otimes\Hom_A(L,A)\to J$ is an isomorphism, since
it is a surjection between two isomorphic finitely generated projective
$k$-modules. We have constructed $L=\beta(J,P)\in M(\CC)$
and we have proved that $\beta\alpha=\id$.
Since $\Hom_A(L,A)\otimes_B P\simeq \Hom_k(P,k)\otimes_B P\simeq k$,
it follows that $\alpha\beta=\id$.

The last assertion of the Proposition is an immediate consequence
of Lemma \ref{orthogonal}.
\end{proof}

\begin{rem}
\label{Serre}
From the previous Theorem, we see that 
$(A/\im\tau_L)\mMod$ is a Serre subcategory of $A\mMod$
(\ie, closed under extensions, subobjects and quotients).
\end{rem}

Let us now study the relation between projective $A$-modules and
projective $(A/J)$-modules.

\begin{lemma}
\label{quoproj}
Let $L\in M(\CC)$ and $J=\im\tau_L$.

Given $P\in\CC\mproj$, then
$\im\tau'_{L,P}=JP$ and $P/JP$ is a projective $A/J$-module.

Conversely, let $Q\in (A/J)\mproj$. Let $U\in k\mproj$ and
$f:U\to \Ext^1_A(Q,L)$ be a surjection. Let $0\to L\otimes U^*\to P\to Q\to 0$
be the extension corresponding to $f$ via the canonical isomorphism
$\Hom_k(U,\Ext^1_A(Q,L))\iso \Ext^1_A(Q,L\otimes U^*)$. Then, $P\in\CC\mproj$.
\end{lemma}

\begin{proof}
The first assertion reduces to the case $P=A$, where it is clear.

Let us now consider the second assertion.
It reduces to the case $Q=(A/J)^n$ for some positive integer $n$.
The canonical map $A^n\to (A/J)^n$ factors through
$\phi:A^n\to P$.
Let $\psi=\phi+\can:A^n\oplus L\otimes U^*\to P$ and $N=\ker\psi$.
Then, $\psi$ is surjective and there is an exact sequence of $A$-modules
$0\to N\to J^{\oplus n}\oplus L\otimes U^*\to L\otimes U^*\to 0$.
Such a sequence splits, hence $N\simeq L\otimes V$ for some $V\in k\mproj$.
By construction, $\Ext^1_A(P,L)=0$, so $\Ext^1_A(P,N)=0$. It follows
that $\psi$ is a split surjection, hence $P$ is projective.
\end{proof}

The following Lemma shows that $M(\CC)$ behaves well with respect to base
change.
\begin{lemma}
\label{basechange}
Let $L$ be an object of $\CC$.
Let $R$ be a commutative noetherian $k$-algebra. If $L\in M(\CC)$,
then $RL\in M(R\CC)$.

The following assertions are equivalent
\begin{itemize}
\item[(i)] $L\in M(\CC)$.
\item[(ii)] $k_\Gm L\in M(k_\Gm\CC)$ for every maximal
ideal $\Gm$ of $k$.
\item[(iii)] $L$ is projective over $k$ and $L(\Gm)\in M(\CC(\Gm))$ 
for every maximal ideal $\Gm$ of $k$.
\end{itemize}
\end{lemma}

\begin{proof}
There is a commutative diagram
$$\xymatrix{
RL\otimes_R\Hom_{RA}(RL,RA)\ar[rr]^-{\tau'_{RL,RA}}&&
RA \\
R(L\otimes \Hom_A(L,A))\ar[urr]_{R\tau'_{L,A}}\ar[u]^\sim
}$$
This shows the first assertion.

Assume (ii). Since $k_\Gm L$ is a projective
$k_\Gm A$-module for every $\Gm$, it follows that $L$ is projective $A$-module.
We obtain also that $\tau'_{L,A}$ is injective and that its cokernel
is projective over $k$. So, (ii)$\Longrightarrow$(i).

Assume (iii). Then, $L$ is a projective $A$-module. Also, 
$\tau'_{L,A}$ is injective. The exact sequence
$0\to L\otimes \Hom_A(L,A)\to L\to \coker \tau'_{L,A}\to 0$ remains exact
after tensoring by $k/\Gm$ for every $\Gm$, hence
$\Tor_1^k(k/\Gm,\coker \tau'_{L,A})=0$ for all $\Gm$, so
$\coker \tau'_{L,A}$ is projective over $k$. Hence, (iii)$\Longrightarrow$(i).

Finally, (i)$\Longrightarrow$(ii) and (i)$\Longrightarrow$(iii) are special
cases of the first part of the Lemma.
\end{proof}

\subsubsection{Definition}
Let $\CC$ be (a category equivalent to)
the module category of a finite projective $k$-algebra $A$.
Let $\Delta$ be a finite set of objects of $\CC$ together with a poset
structure.

Given $\Gamma$ an ideal of $\Delta$, we
denote by $\CC[\Gamma]$ the full subcategory of $\CC$ of objects 
$M$ such that $\Hom(D,M)=0$ for all $D\in \Delta\setminus\Gamma$.

We put $\tDelta=\{D\otimes U|D\in\Delta,\ U\in k\mproj\}$. We put the order
on $\tDelta$ given by $D\otimes U<D'\otimes U'$ if $D<D'$.
We also put $\Delta_\otimes=\{D\otimes U|D\in\Delta,\ U\in \Pic(k)\}$.

\begin{defi}
\label{defhwcat}
We say that $(\CC,\Delta)$ is a highest weight
category over $k$ if the following conditions are satisfied:
\begin{itemize}
\item[(1)]
The objects of $\Delta$ are projective over $k$.
\item[(2)]
$\End_\CC(M)=k$ for all $M\in\Delta$.
\item[(3)]
Given $D_1,D_2\in\Delta$, if
$\Hom_\CC(D_1,D_2)\not=0$, then $D_1\le D_2$.
\item[(4)]
$\CC[\emptyset]=0$.
\item[(5)]
Given $D\in\Delta$, there is $P\in\CC\mproj$ and $f:P\to D$ surjective
such that $\ker f\in \CC^{\tDelta_{>D}}$.
\end{itemize}
\end{defi}

We call $\Delta$ the set of {\em standard objects}.

\smallskip
Let $(\CC,\Delta)$ and $(\CC',\Delta')$ be two highest weight categories
over $k$. A functor $F:\CC\to\CC'$ is an {\em equivalence of highest weight
categories} if it is an equivalence of categories and
if there is a bijection $\phi:\Delta\iso\Delta'$ and 
invertible $k$-modules $U_D$ for $D\in\Delta$ such that
$F(D)\simeq \phi(D)\otimes U_D$ for $D\in\Delta$.

When $k$ is a field, this corresponds to the usual concept of a highest
category \cite{CPS1}. We leave it to the interested reader to extend the
definition to the case where $\Delta$ is an infinite set (this will cover
representations of reductive groups over $\BZ$) and to the non split situation
where (2) is relaxed.

\begin{lemma}
\label{inductive}
Let $\CC$ be the module category of a finite projective $k$-algebra.
Let $\Delta$ be a finite set of objects of $\CC$ together with a poset
structure. Let $L$ be a maximal element of $\Delta$.

Then,
$(\CC,\Delta)$ is a highest weight category if and only if
$L\in M(\CC)$ and $(\CC[\Delta\setminus\{L\}],\Delta\setminus\{L\})$ is
a highest weight category.
\end{lemma}

\begin{proof}
Assume $(\CC,\Delta)$ is a highest weight category.
Given $D\in\Delta$, let $P_D$ be a projective object of $\CC$ with
a surjection $P_D\to D$ whose kernel is in $\CC^{\tDelta_{>D}}$
(Definition \ref{defhwcat} (5)). Let $P=\bigoplus_{D\in\Delta}P_{D}$.
Then, $P$ is a progenerator for $\CC$ (Definition \ref{defhwcat} (4)).

By Definition \ref{defhwcat} (5), $L$ is projective.
We deduce that $P_{D}$ has a submodule
$Q_{D}\simeq L\otimes U_{D}$ for some $U_{D}\in k\mproj$ with
$P_{D}/Q_{D}\in \CC^{(\widetilde{\Delta\setminus\{L\}})_{>D}}$. So,
$P$ has a submodule $Q\simeq L\otimes U$ for some $U\in k\mproj$
with $P/Q\in \CC^{\widetilde{\Delta\setminus\{L\}}}\subset
 \CC[\Delta\setminus\{L\}]$.
We deduce that $L\in M(\CC)$. Also, 
$P_{D}/Q_{D}$ is a projective object of $\CC[\Delta\setminus\{L\}]$
(Lemma \ref{quoproj})
and (5) holds for $\CC[\Delta\setminus\{L\}]$.
So, $(\CC[\Delta\setminus\{L\}],\Delta\setminus\{L\})$ is
a highest weight category.

\smallskip
Assume now $L\in M(\CC)$ and
$(\CC[\Delta\setminus\{L\}],\Delta\setminus\{L\})$ is
a highest weight category.

Let $D\in\Delta\setminus\{L\}$ and $Q$ be a projective object of
$\CC[\Delta\setminus\{L\}]$ as in Definition \ref{defhwcat} (5).
Let $U\in k\mproj$ and $p:U\to \Ext^1_\CC(Q,L)$ be a surjection.
Via the canonical isomorphism $\Hom_k(U,\Ext^1_\CC(Q,L))\iso
\Ext^1_\CC(Q,L\otimes U^*)$, this gives an extension
$0\to L\otimes U^*\to P\to Q\to 0$. By Lemma \ref{quoproj}, 
$P$ is projective (in $\CC$). So, (5) holds for $\CC$ and
$\CC$ is a highest weight category.
\end{proof}

\begin{prop}
\label{prophwcat}
Let $(\CC,\Delta)$ be a highest weight category.
Then,
\begin{itemize}
\item
Given $\Gamma$ an ideal of $\Delta$, then $(\CC[\Gamma],\Gamma)$ is
a highest weight category and $\CC[\Gamma]$ is the full subcategory of
$\CC$ with objects the quotients of objects of $\CC^{\tGamma}$. This
is a Serre subcategory of $\CC$.

\item
Given $D_1,D_2\in\Delta$, if
$\Ext^i_\CC(D_1,D_2)\not=0$ for some $i$, then $D_1\le D_2$.
Furthermore, $\Ext^i_\CC(D_1,D_1)=0$ for $i>0$.
\item
Let $P\in\CC\mproj$ and let
$\Delta\iso\{1,\ldots,n\},\ \Delta_i\leftrightarrow i$, be an
increasing bijection. Then, there is
a filtration $0=P_{n+1}\subset P_{n}\subset\cdots
\subset P_1=P$ with $P_i/P_{i+1}\simeq \Delta_i\otimes U_i$ for
some $U_i\in k\mproj$.
\end{itemize}
\end{prop}

\begin{proof}
By induction, it is sufficient to prove the first assertion in
the case where $|\Delta\setminus\Gamma|=1$. It is then given by Lemma
\ref{inductive}.

Let us now prove the second assertion. Let $\Omega$ be a coideal of
$\Delta$. Then, every object of $\CC^{\tOmega}$ has a projective resolution
with terms in $\CC^{\tOmega}$. This shows the first part of the second
assertion.
The second part follows from the fact that there is a projective $P$
and $f:P\to D_1$ surjective with kernel in $\CC^\tOmega$, with
$\Omega=\Delta_{>D_1}$.

The last assertion follows easily by induction on $|\Delta|$
from Lemma \ref{inductive} and its proof.
\end{proof}

\begin{prop}
\label{hwbasechange}
Let $k'$ be a commutative noetherian $k$-algebra. Let $(\CC,\Delta)$ be
a highest weight category over $k$. Then $(k'\CC,k'\Delta)$ is
a highest weight category over $k'$ and
$(k'\CC)[k'\Gamma]\simeq k'(\CC[\Gamma])$ for all ideals $\Gamma$
of $\Delta$.
\end{prop}

\begin{proof}
Let $A$ be a finite projective $k$-algebra with an equivalence
$\CC\iso A\mMod$.
Let $L$ be a maximal element of $\Delta$. Then, $L\in M(\CC)$
(Lemma \ref{inductive}) and $k'L\in M(k'\CC)$
(Lemma \ref{basechange}). Let $J$ be the ideal of $A$ corresponding to
$L$. For $\Gamma=\Delta\setminus\{L\}$, we have $\CC[\Gamma]\iso (A/J)\mMod$,
$(k'\CC)[k'\Gamma]\iso k'(A/J)\mMod$,
and we deduce that $(k'\CC)[k'\Gamma]\simeq k'(\CC[\Gamma])$.

The Proposition follows by induction on $|\Delta|$ from Lemmas \ref{basechange}
and \ref{inductive}.
\end{proof}

Testing that $(\CC,\Delta)$ is a highest weight category can be reduced to
the case of a base field:
\begin{thm}
\label{carhwcat}
Let $\CC$ be the module category of a finite projective $k$-algebra.
Let $\Delta$ be a finite poset of objects of $\CC\cap k\mproj$.

Then, $(\CC,\Delta)$ is a highest weight category if and only
if $(\CC(\Gm),\Delta(\Gm))$ is a highest weight category
for all maximal ideals $\Gm$ of $k$.
\end{thm}

\begin{proof}
The first implication is a special case of Proposition \ref{hwbasechange}.
The reverse implication follows by induction on $|\Delta|$ from Lemmas
\ref{basechange} and \ref{inductive}.
\end{proof}

\subsubsection{Quasi-hereditary algebras}
Let us recall now the definition of split quasi-hereditary algebras
\cite[Definition 3.2]{CPS2}.

A structure of {\em split quasi-hereditary algebra} on a finite projective
$k$-algebra $A$
is the data of a poset $\Lambda$ and of a set of ideals
$\CI=\{I_\Omega\}_{\Omega\text{ coideal of }\Lambda}$ of $A$ such that
\begin{itemize}
\item given $\Omega\subset\Omega'$ coideals of $\Lambda$, then
$I_\Omega\subset I_{\Omega'}$
\item given $\Omega\subset\Omega'$ coideals of $\Lambda$ with
$|\Omega'\setminus\Omega|=1$, then 
$I_{\Omega'}/I_\Omega$ is an indecomposable split heredity ideal of
$A/I_\Omega$
\item $I_\emptyset=0$ and $I_\Lambda=A$.
\end{itemize}

The following Theorem shows that notion of highest weight category 
corresponds to that of split quasi-hereditary algebras.
\begin{thm}
\label{hwqh}
Let $A$ be a finite projective $k$-algebra and let $\CC=A\mMod$.

Assume $A$, together with $\Lambda$ and $\CI$ is a split quasi-hereditary
algebra. Given $\lambda\in\Lambda$, let
$\Delta(\lambda)\in M((A/I_{\Lambda_{>\lambda}})\mMod)$
correspond to $I_{\Lambda_{\ge \lambda}}/I_{\Lambda_{>\lambda}}$.
Then, $(\CC,\{\Delta(\lambda)\}_{\lambda\in\Lambda})$
is a highest weight category.

Conversely, assume $(\CC,\Delta)$ is a highest weight
category. Given $\Omega$ a coideal of $\Delta$, let
$I_\Omega\subset A$ be the annihilator of all objects of 
$\CC[\Delta\setminus\Omega]$.
Then, $A$ together with $\{I_\Omega\}_\Omega$ is a split quasi-hereditary
algebra and $(A/I_\Omega)\mMod$ identifies with $\CC[\Delta\setminus\Omega]$.
\end{thm}

\begin{proof}
We prove the first assertion by induction on $|\Lambda|$.
Assume $A$ is a split quasi-hereditary algebra.
Let $\lambda\in\Lambda$ be maximal and let $\Gamma=\Lambda\setminus\{\lambda\}$.
 Let $J=I_{\lambda}$. By Proposition \ref{moduleideal},
we have $\CC[\{\Delta(\lambda')\}_{\lambda'\in\Gamma}]\iso (A/J)\mMod$. Since $A/J$ is a
split quasi-hereditary
algebra, it follows by induction that
$(\CC[\{\Delta(\lambda')\}_{\lambda'\in\Gamma}],\{\Delta(\lambda')\}_{\lambda'\in\Gamma})$
is a highest weight
category. By Lemma \ref{inductive}, it follows that
$(\CC,\{\Delta(\beta)\}_{\beta\in\Lambda})$ is a highest weight category.

\smallskip
We prove the second assertion by induction on $|\Delta|$.
Let $(\CC,\Delta)$ be a highest weight category.
Let $\Omega\subset\Omega'$ be coideals of $\Delta$ with $|\Omega'\setminus
\Omega|=1$. If $\Omega=\emptyset$, then $\Omega'=\{L\}$ and $L\in M(\CC)$,
hence $I_{\{L\}}$ is an indecomposable split heredity ideal of $A$
(Proposition \ref{moduleideal}).
Assume now $\Omega\not=\emptyset$ and let $L$ be a maximal
element of $\Omega$. Then, $\CC[\Delta\setminus\{L\}]\simeq (A/I_{\{L\}})\mMod$
(Proposition \ref{moduleideal}). By induction, 
$I_{\Omega'}/I_\Omega$ is an indecomposable split heredity ideal of
$A/I_\Omega$. So, $A$ is a split quasi-hereditary algebra.
\end{proof}

\begin{rem}
Note that, starting from a split quasi-hereditary algebra, we obtain a well
defined poset $\Delta_\otimes$, but
$\Delta$ is not unique, unless $\Pic(k)=1$.
\end{rem}

\begin{rem}
Note that Theorem \ref{carhwcat} translates, via Theorem \ref{hwqh}, to
a known criterion for split quasi-heredity \cite[Theorem 3.3]{CPS2}.
\end{rem}

\subsubsection{Tilting objects}

\begin{prop}
\label{defnabla}
Let $(\CC,\{\Delta(\lambda)\}_{\lambda\in\Lambda})$ be a highest weight
category. Then, there
is a set $\{\nabla(\lambda)\}_{\lambda\in\Lambda}$
of objects of $\CC$, unique up to isomorphism, with the following
properties
\begin{itemize}
\item
$(\CC^\opp,\{\nabla(\lambda)\}_{\lambda\in\Lambda})$ is a highest weight
category.
\item
Given $\lambda,\beta\in\Lambda$, then
$\Ext^i_\CC(\Delta(\lambda),\nabla(\beta))\simeq
\begin{cases}
k & \text{ if }i=0 \text{ and }\lambda=\beta\\
0 & \text{ otherwise.}
\end{cases}$
\end{itemize}
\end{prop}

\begin{proof}
Let $A$ be a finite projective $k$-algebra with $A\mMod\iso\CC$, together
with its structure $\CI$ of split quasi-hereditary algebra (Theorem \ref{hwqh}).
Then, $A^\opp$ together with $\CI$ is a split quasi-hereditary algebra
\cite[Corollary 3.4]{CPS2}. Let $\CC^*=A^\opp\mMod$ and
$\{\Delta(\lambda^*)\}_{\lambda\in\Lambda}$ be a corresponding set
of standard objects.

We have $\Ext^{>0}_A(\Delta(\lambda),\Delta(\beta^*)^*)=0$
for all $\beta$ with $\beta\not>\lambda$, since 
$\Delta(\lambda),\Delta(\beta^*)^*\in (A/I_{>\lambda})\mMod$
and $\Delta(\lambda)$ is a projective $(A/I_{>\lambda})$-module. Similarly,
we have $\Ext^{>0}_{A^\opp}(\Delta(\beta^*),\Delta(\lambda)^*)=0$
if $\lambda\not>\beta$. Since
$\Ext^{>0}_A(\Delta(\lambda),\Delta(\beta^*)^*)^*\simeq
\Ext^{>0}_{A^\opp}(\Delta(\beta^*),\Delta(\lambda)^*)$, we deduce
that this vanishes for all $\lambda,\beta$.
In the same way, we obtain
$\Hom_A(\Delta(\lambda),\Delta(\beta^*)^*)=0$ for $\beta\not=\lambda$.

Let $\Gm$ be a maximal ideal of $k$. We know that
$\Hom_{A(\Gm)}(\Delta(\lambda)(\Gm),\Delta(\lambda^*)(\Gm)^*)=
k/\Gm$ (cf eg \cite[proof of Theorem 3.11]{CPS1}). Let
$U_\lambda=\Hom_A(\Delta(\lambda),\Delta(\lambda^*)^*)$. This is
a projective $k$-module, since
$\Delta(\lambda),\Delta(\lambda^*)^*\in (A/I_{>\lambda})\mMod$,
$\Delta(\lambda^*)^*$ is a projective $k$-module,
and $\Delta(\lambda)$ is a projective $(A/I_{>\lambda})$-module.
It follows that $U_\lambda$ is invertible.
Let $\nabla(\lambda)=U_\lambda^*\otimes\Delta(\lambda^*)^*$.
Then, $\Hom_A(\Delta(\lambda),\nabla(\lambda))\simeq k$.

\smallskip
Let us now show the unicity part. Let
$\{\nabla'(\lambda)\}_{\lambda\in\Lambda}$ be a set of objects of $\CC$
with the same properties. We show by induction that
$\nabla'(\lambda)\simeq\nabla(\lambda)$.

Assume this holds for $\lambda>\alpha$. Then,
$\{\nabla'(\lambda)^*\}_{\lambda\not>\alpha}$ and
$\{\nabla(\lambda)^*\}_{\lambda\not>\alpha}$ are sets of standard objects
for a highest weight category structure on $(A/I_{>\alpha})^\opp\mMod$.
The maximality of $\alpha$ shows that $\nabla'(\alpha)^*$ is a
projective $(A/I_{>\alpha})^\opp$-module, hence it has a filtration
$0=P_n\subset\cdots\subset P_1=\nabla'(\alpha)^*$ such that
$P_i/P_{i+1}\simeq \nabla(\lambda_i)^* \otimes U_i$ for some
$U_i\in k\mproj$ and $\lambda_i{\not>}\alpha$,
as in Proposition \ref{prophwcat}. By assumption, we have
$\Hom(\nabla'(\alpha)^*,\Delta(\beta)^*\otimes U)\simeq
\Hom(\Delta(\beta),\nabla'(\alpha))\otimes U=\delta_{\alpha,\beta}\cdot U$ and
$\Ext^1(\nabla'(\alpha)^*,\Delta(\beta)^*\otimes U)=0$, hence
there is a unique term in the filtration and
$\nabla'(\alpha)\simeq\nabla(\alpha)$.
\end{proof}

We put $\nabla=\{\nabla(\lambda)\}_{\lambda\in\Lambda}$ and
$\tnabla=\{L\otimes U|L\in\nabla,\ U\in k\mproj\}$.

\smallskip
From Proposition \ref{defnabla} and its proof, we deduce:
\begin{prop}
Given $\lambda\in\Lambda$, there is a relatively $k$-injective
module $I$ and an injection $g:\nabla(\lambda)\to I$ with
$\coker g\in \CC^{\tnabla_{>\lambda}}$.
\end{prop}

\begin{lemma}
\label{carfilt}
Let $M\in\CC\cap k\mproj$. Then, $M\in\CC^{\tDelta}$ if and only
if $\Ext^1_\CC(M,\nabla(\lambda))=0$ for all $\lambda\in\Lambda$.
Similarly, $M\in\CC^{\tnabla}$ if and only
if $\Ext^1_\CC(\Delta(\lambda),M)=0$ for all $\lambda\in\Lambda$.
\end{lemma}

\begin{proof}
The first implication is clear. We show the converse by induction on
$|\Lambda|$.
Let $M\in\CC\cap k\mproj$ with $\Ext^1_\CC(M,\nabla(\lambda))=0$ for all
$\lambda\in\Lambda$

Let $\lambda\in\Lambda$ be maximal. Let $M_0=\im\tau_{\Delta(\lambda),M}$,
a subobject of $M$ together with a surjective map
$f:\Delta(\lambda)\otimes U\to M_0$, where $U=\Hom_\CC(\Delta(\lambda),M)
\in k\mproj$. Given $\lambda'\not=\lambda$, we have
$\Hom_\CC(M_0,\nabla(\lambda'))=0$, hence
$\Ext^1_\CC(M/M_0,\nabla(\lambda'))=0$.
We have $M/M_0\in \CC[\{\Delta(\lambda')\}_{\lambda'\not=\lambda}]$, hence
$M/M_0\in\CC^{\tDelta}$ by induction. So,
$\Ext^i_\CC(M/M_0,\nabla(\lambda'))=0$ for all $i>0$ and $\lambda'\in\Lambda$.
Consequently, 
$\Ext^1_\CC(M_0,\nabla(\lambda'))=0$ for all $\lambda'\in\Lambda$.

Let $N=\ker f$. We have $\Hom_\CC(\Delta(\lambda),\nabla(\lambda'))=0$
for $\lambda'\not=\lambda$, hence $\Hom_\CC(N,\nabla(\lambda'))=0$
for $\lambda'\not=\lambda$.

By construction, the canonical map
$\Hom_\CC(\Delta(\lambda),\Delta(\lambda)\otimes U)\to
 \Hom_\CC(\Delta(\lambda),M_0)$
is surjective. So, $\Hom_\CC(\Delta(\lambda),N)=0$. Let $P$
be a projective object of $\CC$ with a surjection $P\to N$. There is 
a subobject $P_0$ of $P$ with $P_0\simeq \Delta(\lambda)\otimes U'$ for some
$U'\in k\mproj$ and $P/P_0 \in
 \CC^{\widetilde{\{\Delta(\lambda')\}}_{\lambda'\not=
\lambda}}$. We obtain a surjection $P/P_0\to N$. We have
$\Hom_\CC(P/P_0,\nabla(\lambda))=0$, hence $\Hom_\CC(N,\nabla(\lambda))=0$.
We deduce that $N=0$, hence $M\in\CC^{\tDelta}$.

The second statement follows by duality.
\end{proof}

The following Lemma will be useful in \S \ref{secunicity}.

\begin{lemma}
\label{charproj}
Let $A$ be a split quasi-hereditary $k$-algebra. Let $M\in (A\mMod)^\tDelta$.
If $\Ext^1_A(M,N)=0$ for all $N\in\Delta$, then $M$ is projective.
\end{lemma}

\begin{proof}
We have $\Ext^1_A(M,N)=0$ for all $N\in (A\mMod)^\tDelta$.
Let $0\to N\to P\to M\to 0$ be an exact sequence with $P$ projective.
Then, $N$ is $\tDelta$-filtered (Lemma \ref{carfilt}), hence
$\Ext^1_A(M,N)=0$ and the sequence splits.
\end{proof}

Recall that the category of perfect complexes for $A$ is the full
subcategory of $D^b(A\mMod)$ of objects isomorphic to
a bounded complex of finitely generated projective $A$-modules.

\begin{prop}
\label{finiteprojdim}
Every object of $\CC\cap k\mproj$ has finite projective dimension.
More precisely, a complex of $\CC$ that is perfect as a complex of
$k$-modules is also perfect as a complex of $\CC$.
\end{prop}

\begin{proof}
This is almost \cite[Theorem 3.6]{CPS2}, whose proof we follow.
We show the Proposition by induction on $|\Lambda|$. Consider
$\lambda\in\Lambda$
maximal and let $J$ be the ideal of $A$ corresponding to the
projective object $L=\Delta(\lambda)$.
Note that we have an isomorphism of $(A,A)$-bimodules
$L\otimes L^\vee\iso J$, where $L^\vee=\Hom_A(L,A)$.
The exact sequence of $(A,A)$-bimodules
$$0\to J\to A\to A/J\to 0$$
induces an exact sequence of functors $A\mMod\to A\mMod$
$$0\to L\otimes \Hom_A(L,-)\to \Id\to (A/J)\otimes_A -\to 0.$$
Let $C$ be a complex of $A$-modules.
We have a distinguished triangle
$$L\otimes \Hom_A(L,C)\to C\to A/J\otimes^{\BL}_A C\rightsquigarrow.$$
Assume $C$ is perfect, viewed as a complex of $k$-modules.
Then, $\Hom_A(L,C)$ is perfect as a complex of $k$-modules, hence
$L\otimes \Hom_A(L,C)$ is perfect as a complex of $A$-modules.
In particular, $A/J\otimes^{\BL}_A C$ is an object of
$D^b((A/J)\mMod)$ that is perfect as a complex of $k$-modules. By induction,
it is perfect as a complex of $(A/J)$-modules. Since
$A/J$ is perfect as a complex of $A$-modules, it follows that
$A/J\otimes^{\BL}_A C$ is a perfect complex of $A$-modules, hence $C$
as well.
\end{proof}

\begin{rem}
Let $\CT$ be the full subcategory of $D^b(\CC)$ of complexes that are
perfect as complexes of $k$-modules. 
Fix an increasing bijection $\Lambda\iso \{1,\ldots,n\}$.
Then, $\CT$ has a 
semi-orthogonal decomposition
$\CT=\langle \Delta_1\otimes k\mperf,\Delta_2\otimes k\mperf,\ldots,
\Delta_n\otimes k\mperf\rangle$. This gives an isomorphism
$$K_0(k\mproj)^\Delta\iso K_0(\CT)=K_0(\CC\mproj),\
\{[L_\lambda]\}_{\lambda\in\Lambda}\mapsto 
\sum_\lambda [\Delta(\lambda)\otimes L]$$
In the isomorphism above, one can replace $\Delta(\lambda)$ by a projective
object $P(\lambda)$ as in Definition \ref{defhwcat}, (5), or
by $\nabla(\lambda)$, $I(\lambda)$ or $T(\lambda)$.
We recover \cite[Corollary 1.2.g]{Do2} (case of integral Schur algebras).
\end{rem}

\begin{defi}
An object $T\in\CC$ is {\em tilting} if $T\in\CC^{\tDelta}\cap \CC^{\tnabla}$.
We denote by $\CC\mtilt$ the full subcategory of $\CC$ of tilting objects.
\end{defi}

\smallskip
A tilting complex is a perfect complex $C$ with the following properties
\begin{itemize}
\item
$C$ generates the category of perfect complexes as a full triangulated
subcategory closed under taking direct summands and
\item
$\Hom_{D^b(A\mMod)}(C,C[i])=0$ for $i\not=0$.
\end{itemize}
Note that a tilting module
is {\em not} a tilting complex in general, for the generating property will be
missing in general. Nevertheless, there is a tilting module which is a tilting
complex, as explained below.

\begin{prop}
\label{tilting}
Let $M\in\CC^{\tDelta}$. Then, there is $T\in\CC\mtilt$
and an injection $i:M\to T$ with $\coker i\in \CC^{\tDelta}$.

Let $\lambda\in\Lambda$. 
There is $T(\lambda)\in\CC\mtilt$ and
\begin{itemize}
\item an injection
$i:\Delta(\lambda)\to T(\lambda)$ with $\coker i\in\CC^{\tDelta_{<\lambda}}$;
\item
a surjection
$p:T(\lambda)\to \nabla(\lambda)$ with $\ker p\in\CC^{\tnabla_{<\lambda}}$.
\end{itemize}
Let $T=\bigoplus_{\lambda\in\Lambda}T(\lambda)$.
Then, $T$ is a tilting complex. Let $A^r=\End_\CC(T)$ and $\CC^r=A^r\mMod$.
There is an equivalence $\BR\Hom_\CC(T,-):D^b(\CC)\iso D^b(\CC^r)$.
Let $\Lambda^r=\{\lambda^r\}_{\lambda\in\Lambda}$ be the opposite
poset to $\Lambda$.
Let $\Delta(\lambda^r)=\BR\Hom_\CC(T,\nabla(\lambda))$. Then,
$(\CC^r,\{\Delta(\lambda^r)\}_{\lambda^r\in\Lambda^r})$ is a
highest weight category.
\end{prop}

\begin{proof}
Let us fix an increasing bijection
$\Delta\iso\{1,\ldots,n\},\ \Delta_i\leftrightarrow i$.
Let $M\in\CC^{\tDelta}$. We construct by induction an object $T$ with
a filtration $0=T_{n+1}\subset M=T_n\subset\cdots\subset T_0=T$
such that $T_{i-1}/T_i\simeq\Delta_i\otimes U_i$ for some $U_i\in k\mproj$,
for $i\le n$.

Assume $T_i$ is defined for some $i$, $2\le i\le n$.
Let $U_i\in k\mproj$ with a surjection of $k$-modules
$U_i\to \Ext^1_\CC(\Delta_i,T_i)$. Via the canonical isomorphism
$\Hom_k(U_i,\Ext^1_\CC(\Delta_i,T_i))\iso \Ext^1_\CC(\Delta_i\otimes U_i,T_i)$,
we obtain an extension
$$0\to T_i\to T_{i-1}\to \Delta_i\otimes U_i\to 0.$$
By construction, we have $\Ext^1(\Delta_i,T_{i-1})=0$, since
$\Ext^1_\CC(\Delta_i,\Delta_i)=0$ (Proposition \ref{prophwcat}).

We have $T/M,T\in \CC^{\tDelta}$. By Proposition \ref{prophwcat},
we have $\Ext^1_\CC(\Delta_i,T)\simeq \Ext^1_\CC(\Delta_i,T_{i-1})=0$.
It follows from Lemma \ref{carfilt} that $T$ is tilting.

\smallskip
Assume $M=\Delta(\lambda)$. Then, in the construction above, one can
replace the bijection
$\Delta\iso\{1,\ldots,n\},\ \Delta_i\leftrightarrow i$
 by an increasing bijection
$\Delta_{<\lambda}\iso\{1,\ldots,m\},\ \Delta_i\leftrightarrow i$, and
obtain the same conclusion. This produces a tilting object $T(\lambda)$.
It has a $\tnabla$-filtration with top term
$\nabla(\lambda)$ giving rise to a map $p:T(\lambda)\to\nabla(\lambda)$
as required.

Every object of $\CC^{\tDelta}$ has finite homological dimension
(Proposition \ref{finiteprojdim}). In
particular, $T$ is a perfect complex. We have
$\Ext^i_\CC(T,T)=0$ for $i\not=0$ by Lemma \ref{carfilt}.
Let $\CD$ be the smallest full triangulated subcategory of $D^b(\CC)$
containing $T$ and closed under direct summands. By induction, $\CD$
contains $\Delta$, hence it contains the projective objects of $\CC$. So,
$T$ generates the category of perfect complexes and $T$ is a tilting
complex.

As a consequence, we have an equivalence $F=\BR\Hom_\CC(T,-)$. Note
that $\Delta(\lambda^r)$ has non zero homology only in degree $0$ by
Lemma \ref{carfilt} and it is projective over $k$.
Let $P(\lambda^r)=F(T(\lambda))$, an object with
homology concentrated in degree $0$ and projective. Also,
we obtain a surjection $P(\lambda^r)\to\Delta(\lambda^r)$ with
kernel filtered by terms $\Delta(\beta^r)\otimes U$ with
$U\in k\mproj$ and $\beta^r>\lambda^r$. This shows that
$(\CC^r,\{\Delta^r(\lambda)\})$ is a highest weight category.
\end{proof}

The highest weight category $\CC^r$ in the Proposition above is called
the {\em Ringel dual} of $\CC$.

\begin{prop}
Fix a family $\{T(\lambda)\}_{\lambda\in\Lambda}$ as in Proposition
\ref{tilting}. Then, every
tilting object of $\CC$ is a direct summand of a direct sum of
$T(\lambda)$'s.

Furthermore, the category $\CC^r$ is independent of the choice of the
$T(\lambda)$'s, up to equivalence of highest weight categories.
\end{prop}

\begin{proof}
The first assertion follows using these equivalences from the fact that
every projective object of $\CC^r$ is a direct summand of a direct sum of
$P(\lambda^r)$'s.

Consider another family $\{T'(\lambda)\}$ and the
associated $T',\CC^{\prime r}$.
We consider the composite equivalence
$$F:\BR\Hom_\CC(T',-)\circ \BR\Hom_\CC(T,-)^{-1}:
D^b(\CC^r)\iso D^b(\CC^{\prime r}).$$
It sends $\Delta(\lambda^r)$ to $\Delta(\lambda^{\prime r})$, hence
it sends projective objects to objects with homology only in degree $0$,
which are projective by Lemma \ref{carfilt}. So, $F$ restricts
to an equivalence of highest weight categories $\CC^r\iso \CC^{\prime r}$.
\end{proof}

\begin{rem}
Note that we don't construct a canonical $T(\lambda)$ (nor a canonical
$P(\lambda)$), our construction depends on the choice of
projective $k$-modules mapping onto certain $\Ext^1$'s.
\end{rem}

\begin{rem}
The theory of tilting modules has been developed by Donkin for algebraic groups
over $\BZ$, cf \cite[Remark 1.7]{Do1}.
\end{rem}

\subsubsection{Reduction to fields}

\begin{prop}
\label{testfield}
Let $(\CC,\Delta)$ be a highest weight category over $k$.
Let $M\in\CC\cap k\mproj$. Then,
$M\in \CC^{\tDelta}$ (resp. $M\in\CC\mtilt$) if and only
if the corresponding property holds for $M(\Gm)$ in
$\CC(\Gm)$, for all maximal ideals $\Gm$ of $k$.

If $k$ is indecomposable, then the same statement holds for the properties
of belonging to $\tDelta$ or $\Delta_\otimes$.
\end{prop}

\begin{proof}
Given $V$ a $k$-module, we put $\bar{V}=(k/\Gm)\otimes V$.
Let $C$ be a bounded complex of projective objects of $\CC$ and
$N\in \CC\cap k\mproj$.
Then, we have a canonical isomorphism
$$(k/\Gm)\Hom_{D^b(A)}(C,N)\iso
\Hom_{D^b(\bar{A})}(\bar{C},\bar{N})$$
(this only needs to be checked for $C=A[i]$, where is it clear).
It follows from Proposition \ref{finiteprojdim} that
we have a canonical isomorphism
$$(k/\Gm)\Hom_A(L,N)\iso \Hom_{\bar{A}}(\bar{L},\bar{N})$$
for $L,N\in \CC\cap k\mproj$.

\smallskip
Let $M\in \CC\cap k\mproj$ with $M(\Gm)\in \CC(\Gm)^{\tDelta(\Gm)}$
for every maximal ideal $\Gm$.
We show that $M\in \CC^{\tDelta}$ by induction on the projective dimension 
of $M$ (which is finite by Proposition \ref{finiteprojdim}).
Let
$0\to L\to P\to M\to 0$ be an exact sequence with $P$ projective.
By Lemma \ref{carfilt}, $L(\Gm)\in \CC(\Gm)^{\tDelta(\Gm)}$ for every $\Gm$.

By induction, it follows that $L\in \CC^{\tDelta}$.
Let $N\in\nabla$.
We have $\Ext^{>0}_{A(\Gm)}(L(\Gm),N(\Gm))=0$.
Let $0\to C^r\to\cdots\to C^0\to L\to 0$ be a projective resolution.
Let $D=0\to \Hom_A(C^0,N)\to\Hom_A(C^1,N)\to\cdots\to\Hom_A(C^r,N)\to 0$.
We have $H^i(D(\Gm))\iso \Ext^i_{A(\Gm)}(L(\Gm),N(\Gm))=0$ for
$i>0$. It follows that the complex $D$ is homotopy equivalent to
$H^0(D)$, as a complex of $k$-modules, and $H^0(D)$ is projective.
So, the canonical map
$\Hom_A(L,N)(\Gm)\to \Hom_{A(\Gm)}(L(\Gm),N(\Gm))$ is an isomorphism.

We have a commutative diagram whose horizontal sequences
are exact
$$\xymatrix{
\Hom_A(P,N)(\Gm)\ar[r]\ar[d]^\sim& \Hom_A(L,N)(\Gm)\ar[r]\ar[d]^\sim&
 \Ext^1_A(M,N)(\Gm)\ar[r]\ar[d]&0\\
\Hom_{A(\Gm)}(P(\Gm),N(\Gm))\ar[r]& \Hom_{A(\Gm)}(L(\Gm),N(\Gm))\ar[r]&
 \Ext^1_{A(\Gm)}(M(\Gm),N(\Gm))\ar[r]&0
}$$
We have $\Ext^1_{A(\Gm)}(M(\Gm),N(\Gm))=0$ and it follows that
$\Ext^1_A(M,N)(\Gm)=0$. Since $\Ext^1_A(M,N)$ is a finitely generated
$k$-module, it must thus be $0$. Lemma \ref{carfilt} shows that
$M\in \CC^{\tDelta}$.

The other statements follow easily.
\end{proof}

\begin{rem}
If $k$ is decomposable, then being in $\tDelta$ cannot be tested
locally --- only being a sum of objects of $\tDelta$ can be tested locally.
\end{rem}

\subsection{Covers}
\subsubsection{Double centralizer}
Let $k$ be a commutative noetherian ring and $A$ a finite dimensional $k$-algebra.
Let $\CC=A\mMod$.

Let $P$ be a finitely generated projective $A$-module, $B=\End_A(P)$,
$F=\Hom_A(P,-):A\mMod\to B\mMod$, and
$G=\Hom_B(FA,-):B\mMod\to A\mMod$.
The canonical isomorphism
$\Hom_A(P,A)\otimes_A -\iso \Hom_A(P,-)$ makes $F$ a left adjoint of $G$. We
denote by $\eps:FG\to\Id$ (resp. $\eta:\Id\to GF$)
the corresponding unit (resp. counit). Note that $\eps$ is an isomorphism.

\smallskip
The following Lemma is immediate.
\begin{lemma}
\label{carimG}
Let $M\in A\mMod$. The following assertions are equivalent
\begin{itemize}
\item
the map $\eta(M):M\to GFM$ is an isomorphism
\item
$F$ induces an isomorphism $\Hom_A(A,M)\iso\Hom_B(FA,FM)$
\item $M$ is a direct summand of a module in the image of $G$.
\end{itemize}
\end{lemma}

We will consider gradually stronger conditions on $F$.

\smallskip
Lemma \ref{carimG} gives~:
\begin{prop}
The following assertions are equivalent:
\begin{itemize}
\item the canonical map of algebras $A\to\End_B(FA)$ is an isomorphism
\item for all $M\in A\mproj$, the map $\eta(M):M\to GFM$ is an isomorphism
\item the restriction of $F$ to $A\mproj$ is fully faithful.
\end{itemize}
\end{prop}

Let us name this ``double centralizer'' situation.
\begin{defi}
We say that $(A,P)$ (or $(A\mMod,P)$)
is a {\em cover} of $B$
if the restriction of $\Hom_A(P,-)$ to $A\mproj$ is fully faithful.
We say also that $(\CC,F)$ is a {\em cover} of $B\mMod$.
\end{defi}

\begin{rem}
\label{remMorita}
Let $E=P\otimes_B-:B\mMod\to A\mMod$. This is a left adjoint of $F$.
The canonical map $\Id\to FE$ is an isomorphism.
By Morita theory, the following conditions are equivalent:
\begin{itemize}
\item $F:A\mMod\to B\mMod$ is an equivalence with inverse $G\simeq E$
\item $F:A\mMod\to B\mMod$ is fully faithful
\item for all $M\in A\mMod$, the map $\eta(M):M\to GFM$ is an isomorphism
\item the adjunction map $EFA\to A$ is an isomorphism.
\end{itemize}
\end{rem}

The ``cover'' property can be checked at closed points:
\begin{prop}
\label{pointcover}
Assume $k$ is regular.
If $(A(\Gm),P(\Gm))$ is a cover
of $B(\Gm)$ for every maximal ideal $\Gm$ of $k$, then
$(A,P)$ is a cover of $B$.
\end{prop}

\begin{proof}
Since $(A,P)$ is a cover of $B$ if and only if
$(k_\Gm A,k_\Gm P)$ is a cover of $k_\Gm B$ for every maximal
ideal $\Gm$ of $k$, we can assume $k$ is local.
We prove now the Proposition by induction on the Krull dimension of $k$.
Let $\pi$ be a regular element of the maximal ideal of $k$.
We have a commutative diagram with exact rows
$$\xymatrix{
0\ar[r] & \End_B(FA)\ar[r]^\pi & \End_B(FA)\ar[r] &
 \Hom_B(FA,(k/\pi)FA) \\
 & \End_B(FA)\ar[r]^\pi\ar@{=}[u] & \End_B(FA)\ar[r] \ar@{=}[u] &
 \End_B(FA)\otimes k/\pi\ar[r]\ar@{^{(}->}[u] & 0 \\
0\ar[r] & A\ar[r]^\pi\ar[u] & A\ar[r]\ar[u] & (k/\pi) A\ar[r]\ar[u] & 0
}$$
and the canonical map
$\End_B(FA)\otimes k/\pi\to \Hom_B(FA,(k/\pi)FA)$ is injective.

By induction, $((k/\pi)A,(k/\pi)P)$ is a cover of $(k/\pi)B$, hence the
canonical map $(k/\pi)A\to \Hom_B(FA,(k/\pi)FA)$ is an isomorphism.
It follows that the canonical map $\End_B(FA)\otimes k/\pi\to
\Hom_B(FA,(k/\pi)FA)$ is an isomorphism, hence the canonical map
$(k/\pi)A\to \End_B(FA)\otimes k/\pi$ is an isomorphism as well.
By Nakayama's Lemma, we deduce that the canonical map
$A\to \End_B(FA)$ is an isomorphism.
\end{proof}

\subsubsection{Faithful covers}
We assume now that we are given a highest weight category structure
$(\CC,\Delta)$ on $\CC$. If $\CC$ is a cover of $B\mMod$, we say that
it is a {\em highest weight cover}.

\begin{defi}
Let $i$ be a non-negative integer.
We say that $(A,P)$ (or $(A\mMod,P)$)
is an {\em $i$-faithful cover} of $B$
if $F=\Hom_A(P,-)$ induces isomorphisms
$\Ext^j_A(M,N)\iso \Ext^j_B(FM,FN)$ for all $M,N\in \CC^{\tDelta}$ and $j\le i$.
We say also that $(\CC,F)$ is an {\em $i$-cover} of $B\mMod$.
\end{defi}

\begin{rem}
For $i$ big enough, this will force $F$
to be an equivalence, assuming $k$ is a field.
\end{rem}

\begin{rem}
Note that the $0$-faithfulness assumption is not satisfied in
Soergel's theory on category $\CO$ for a complex semi-simple Lie algebra,
cf already the case of $\Gsl_2$.
\end{rem}

\begin{prop}
\label{charff}
The following assertions are equivalent:
\begin{enumerate}
\item $(\CC,F)$ is a $0$-faithful cover of $B\mMod$.
\item for all $M\in \CC^\tDelta$, the map
$\eta(M):M\to GF M$ is an isomorphism
\item every object of $\CC^\tDelta$ is in the image of $G$
\item for all $T\in \CC\mtilt$, the map
$\eta(T):T\to GF T$ is an isomorphism
\item every object of $\CC\mtilt$ is in the image of $G$.
\end{enumerate}
\end{prop}

\begin{proof}
The equivalence of (1), (2) and (3)
and the equivalence of (4) and (5) is given by Lemma \ref{carimG}.

Assume (4). Let $M\in \CC^\tDelta$. Then there is an exact
sequence
$$0\to M\to T\to N\to 0$$
where $T\in\CC\mtilt$ $N\in\CC^\tDelta$ (Proposition \ref{tilting}).
We have a commutative diagram with exact rows
$$\xymatrix{
0\ar[r] & M \ar[r]\ar[d]^{\eta(M)} & T\ar[r]\ar[d]_\sim^{\eta(T)} &
 N\ar[r]\ar[d]^{\eta(N)} & 0\\
0\ar[r] & GFM \ar[r] & GFT\ar[r] & GFN
}$$
It follows that $\eta(M)$ is injective for all $M\in \CC^\tDelta$.
In particular, in the diagram above, $\eta(N)$ is injective and
it follows that $\eta(M)$ is surjective. So, (4) implies (2)
and the converse is trivial.
\end{proof}

\begin{prop}
\label{charfaithful}
Assume $(\CC,F)$ is a $0$-faithful cover of $B\mMod$.
The following assertions are equivalent:
\begin{enumerate}
\item 
$(\CC,F)$ is a $1$-faithful cover of $B\mMod$.
\item
$F$ restricts to an equivalence of exact categories
$\CC^\tDelta\iso (B\mMod)^{F\tDelta}$ with inverse $G$
\item
for all $M\in \CC^\tDelta$, we have $\BR^1G(FM)=0$
\end{enumerate}
\end{prop}

\begin{proof}
If (2) holds, then $\Ext^1_A(M,N)\iso \Ext^1_B(FM,FN)$ for
$M,N\in (A\mMod)^\tDelta$, \ie, (1) holds.

We have $\BR^1G(FM)=\Ext^1_B(FA,FM)=0$, hence (1)$\Longrightarrow$(3).

Assume (3).
Let $X,Y\in \CC^\tDelta$ and
let $0\to FX\to U\to FY\to 0$ be an exact sequence. We have
an exact sequence $0\to GFX\to GU\to GFY\to 0$. Since 
$X\iso GFX$ and $Y\iso GFY$, we deduce that $GU\in \CC^\tDelta$.
Now, $FGU\iso U$, hence $U\in F(\CC^\tDelta)$. It follows by induction
on the length of a $F\tDelta$-filtration
that $F(\CC^\tDelta)=(B\mMod)^{F\tDelta}$. So, (3) implies (2).
\end{proof}

The following very useful result shows that $1$-faithful quasi-hereditary
covers arise naturally as deformations of $0$-faithful covers.

\begin{prop}
\label{equivDelta}
Assume $k$ is regular and $KA$ is split semi-simple.
If $(A(\Gm),P(\Gm))$ is a $0$-faithful cover of $B(\Gm)$ for
every maximal ideal $\Gm$ of $k$, then
$(A,P)$ is a $1$-faithful cover of $B$.
\end{prop}

\begin{proof}
As in the proof of
Proposition \ref{pointcover}, we can assume $k$ is local with maximal
ideal $\Gm$.

Let us first assume $k$ is a discrete valuation ring with uniformizing
parameter $\pi$.
Let $N\in (A\mMod)^\tDelta$.
The composition of canonical maps
$(k/\pi)N\to (k/\pi)GFN\to G((k/\pi)FN)$ is an isomorphism by
assumption and the second map is surjective, hence
both maps are isomorphisms. By Nakayama's Lemma,
it follows that the canonical map $N\to GFN$ is an isomorphism.
Since $\pi$ is regular for $k$, $FA$ and $FN$,
the Universal Coefficient Theorem (\ie, the isomorphism
$(k/\pi)\otimes_k^\BL \BR G(FN)\iso \BR G(FN\otimes_k^\BL (k/\pi)$)
gives an exact sequence
$$0\to (k/\pi)GFN\to G((k/\pi)FN)\to \Tor_1^k(\BR^1 G(FN),k/\pi)\to 0.$$
We deduce that $\Tor_1^k(E,k/\pi)=0$, where
$E=\BR^1G(FN)$, hence $E$ is free over $k$. 
Note that the canonical map $N'\to GFN'$ is an isomorphism for every
$N'\in KA\mMod$, hence $KB$ is Morita-equivalent to $KA$ (cf Remark
\ref{remMorita}).
Since $KB$ is
semi-simple, $E$ is a torsion $k$-module and this forces $E=0$.
So, the Proposition holds in the case $k$ has Krull dimension $1$.

We prove now the Proposition by induction on the Krull dimension of $k$.
Assume the Krull dimension of $k$ is at least $2$.
There is $\alpha\in k-\{0\}$ such that 
$A[\alpha^{-1}]$ is isomorphic to a product of matrix algebras over
$k[\alpha^{-1}]$. Then
$(k_\Gp/\Gp)A$ is split semi-simple,
whenever $\Gp$ is a prime ideal of $k$ with $\alpha\not\in\Gp$.

We proceed as in the proof of Proposition \ref{pointcover}.
Let $N\in (A\mMod)^{\tDelta}$ such that $\BR^1GFN\not=0$. Let $Z$ be the support
of $\BR^1GFN$ in $\Spec k$, a non-empty strict closed subvariety.
Let $\pi\in \Gm$ regular with $Z\cap \Spec (k/\pi)\not=\emptyset$ and
$\alpha\not\in (\pi)$.

We have a commutative diagram with exact rows
$$\xymatrix{
0\ar[r] & GFN\ar[r]^\pi & GFN\ar[r] &
 G((k/\pi)FN) \\
 & GFN\ar[r]^\pi\ar@{=}[u] & GFN\ar[r] \ar@{=}[u] &
 (k/\pi)GFN\ar[r]\ar@{^{(}->}[u] & 0 \\
0\ar[r] & N\ar[r]^\pi\ar[u] & N\ar[r]\ar[u] &
 (k/\pi) N\ar[r]\ar[u] & 0
}$$
Since the canonical map
$(k/\pi)N\to G((k/\pi)FN)$ is an isomorphism, we deduce
that the canonical map $N\to GFN$ is an isomorphism.
The Universal Coefficient Theorem gives an exact sequence
$$0\to (k/\pi)GFN\to G((k/\pi)FN)\to \Tor_1^k(\BR^1 G(FN),k/\pi)\to 0.$$
It follows that $\BR^1G(FN)$ has no $\pi$-torsion, which is a contradiction.
So, $\BR^1GN=0$. We deduce that
$(A,P)$ is a $1$-faithful cover of $B$.
\end{proof}

\begin{rem}
In the proof above, the case of a discrete valuation ring has been treated
separately, for if $k/\pi$ is finite, then there might be no element
$\alpha$ as needed. An alternative proof would be to take a faithfully flat
extension of $k$ to avoid this problem.
\end{rem}

\subsubsection{Unicity of faithful covers}
\label{secunicity}

\begin{defi}
\label{equivalence}
We say that two highest weight covers $(\CC,F)$ and $(\CC',F')$ of
$B$ are {\em equivalent} if
there is an equivalence of highest weight categories
$\CC\iso \CC'$ making the following diagram commutative
$$\xymatrix{
\CC\ar[dd]_\sim \ar[dr]^-{F} \\
& B\mMod \\
\CC' \ar[ur]_-{F'}
}$$
\end{defi}

The following result shows that a $1$-faithful highest weight cover depends
only on $F\Delta_\otimes$:

\begin{prop}
Let $(\CC,F)$ be a $1$-faithful highest weight cover of $B$.

Given $M\in F\Delta$, there is a pair $(Y(M),p_M)$ unique up to isomorphism
with $Y(M)\in B\mMod$ and
$p_M:Y(M)\to M$ a surjection such that
$\ker p_M\in (B\mMod)^{F\tDelta}$ and
$\Ext^1_B(Y(M),F\Delta)=0$.

Given $N\in\Delta$ with $M=F(N)$ and $q_N:P(N)\to N$ a surjective map
with $\ker q_N\in\CC^{\tDelta}$ and $P(N)$ a projective $A$-module,
then $Y(M)=F(P(N))$ and $p_M=F(q_N)$ satisfy the property above.

\medskip
Let $Y=\bigoplus_{M\in F\Delta}Y(M)$, $A'=\End_B(Y)$,
$\Delta'=\Hom_B(Y,F\Delta)$, and $P'=\Hom_{(A')^\opp}(Y,A')$.
Then, $(A'\mMod,\Delta')$ is a highest weight category and
together with $\Hom_{A'}(P',-)$, this is a
$1$-faithful highest weight cover of $B$ equivalent to $(\CC,F)$.
\end{prop}

\begin{proof}
The unicity follows from Lemma \ref{charproj}, while the construction
of $(Y(M),p_M)$ with the required properties is clear.

Note that $\bigoplus_{N\in\Delta}P(N)$ is a progenerator for $A$, since
every object of $\Delta$ appears as a quotient.
We have a canonical isomorphism $\End_A(\bigoplus_{N\in\Delta}P(N))\iso
A'$, hence an equivalence
$$\Hom_A(\bigoplus_{N\in\Delta}P(N),-): A\mMod\iso A'\mMod$$
 giving rise to the commutative diagram of
Definition \ref{equivalence}.
\end{proof}

We deduce a unicity result.
\begin{cor}
\label{unicityfc}
Let $(\CC,F)$ and $(\CC',F')$ be two $1$-faithful highest weight covers of
$B$. Assume $F\Delta_\otimes\simeq F'\Delta'_\otimes$. Then,
$(\CC,F)$ and $(\CC',F')$ are equivalent highest weight covers.
\end{cor}

\subsubsection{Deformation}
\label{orderqh}
We assume in \S \ref{orderqh} that $k$ is a noetherian domain with field of
fractions $K$.

When $K\CC$ is split semi-simple, we can restate the definition of a highest
weight category structure on $\CC$ as follows (cf \cite[Lemma 1.6]{DuPaSc2}):

\begin{prop}
Let $\CC$ be the module category of a finite projective $k$-algebra
and let $\Delta$ be a finite poset of objects of $\CC\cap k\mproj$.
Assume $K\CC$ is split semi-simple.

Then, $(\CC,\Delta)$ is a highest weight category if and only if
there is a bijection $\Irr(K\CC)\iso\Delta,\ E\mapsto \Delta(E)$, such that
\begin{itemize}
\item $K\Delta(E)\simeq E$ for $E\in\Irr(K\CC)$.
\item for $E\in\Irr(K\CC)$, there is a projective module $P(E)$ 
with a filtration $0=P_r\subset\cdots \subset P_1=P(E)$ such that
$P_1/P_2\simeq \Delta(E)$ and $P_j/P_{j+1}\simeq \Delta(F_j)\otimes U_j$ for
some $F_j>E$ and $U_j\in k\mproj$, for $j\ge 2$.
\item $\bigoplus_{E\in\Irr(K\CC)}P(E)$ is a progenerator of $\CC$.
\end{itemize}
\end{prop}

Note that $\Delta_\otimes$ is determined by the
order on $\Irr(K\CC)$~: given $Q$ a projective 
object of $\CC$ with 
$KQ\simeq E\oplus\bigoplus_{F>E} F^{a_F}$ for some integers
$a_F$, then,
the image of $Q$ by a surjection $KQ\to E$ is isomorphic to
$\Delta(E)\otimes U$ for some $U\in\Pic(k)$.

\medskip

Let $B$ be a finite projective $k$-algebra with $KB$ split semi-simple.
Let $(\CC,F)$ be a $1$-faithful highest weight cover of $B$.
Then, $(K\CC,KF)$ is a $1$-faithful highest weight cover of $KB$, hence
$KF:K\CC\to KB\mMod$ is an equivalence and it induces a bijection
$\Irr(K\CC)\iso\Irr(KB)$. We will say that $(\CC,F)$ is a
highest weight cover of $B$ for the order on $\Irr(KB)$
coming from the one on $\Irr(K\CC)$.
Given $I\subset\Irr(KB)$, we denote by $(KB)^I$ the sum of the simple
$KB$-submodules of $KB$ isomorphic to elements of $I$.

\begin{lemma}
\label{decomfilt}
Let $J\subset I$ be coideals of $\Irr(KB)$
such that no two distinct elements of $I\setminus J$ are comparable.
Then,
$$\left((KB)^I\cap B\right)/\left((KB)^J\cap B\right)\simeq
\bigoplus_{E\in I\setminus J}F\Delta(E)\otimes U_E.$$
where $U_E\in k\mproj$ and $\rank_k U_E=\dim_K E$.
\end{lemma}

\begin{proof}
Recall that $\CC=A\mMod$, $F=\Hom_A(P,-)$ and $B=\End_A(P)$.
Since $P$ is $\tDelta$-filtered, there is a filtration 
$P_0\subset P_1\subset P$ with
$P_0\in \CC^{\tDelta(J)}$,
$P_1/P_0\simeq\bigoplus_{E\in I\setminus J}\Delta(E)\otimes U_E$ for
some $U_E\in k\mproj$, and
$P/P_1\in \CC^{\tDelta(\Irr(KB)\setminus I)}$.
So, we have a filtration $FP_0\subset FP_1\subset FP=B$ and
$(KB)^I\cap B=FP_1$ and $(KB)^J\cap B=FP_0$, since
$FP_0$ and $FP_1$ are direct summands of $FP$ as $k$-modules. Furthermore,
$\dim_K KU_E=\dim_K E$ and we are done.
\end{proof}

We can now show that a $1$-faithful highest weight cover is determined
by the induced order on $\Irr(KB)$.

\begin{thm}
\label{unicityfromorder}
Let $B$ be a finite projective $k$-algebra such that $KB$ is split semi-simple.
Fix two orders, $\le_1$ and $\le_2$ on $\Irr(KB)$.
Let $(\CC_1,F_1)$ and $(\CC_2,F_2)$ be $1$-faithful highest weight covers of
$B$ for the orders $\le_1$ and $\le_2$.

Assume $\le_1$ is a refinement of $\le_2$.
Then, there is an equivalence $\CC_1\mMod\iso \CC_2\mMod$ of highest weight
covers of $B$ inducing the
bijection $\Irr(K\CC_1)\iso\Irr(KB)\iso\Irr(K\CC_2)$.
\end{thm}

\begin{proof}
Let $E\in\Irr(KB)$, $I=\Irr(KB)_{\ge_1 E}$ and
$J=\Irr(KB)_{>_1E}$. These are coideals for $\le_1$ and also for
$\le_2$. Using Lemma \ref{decomfilt}, we obtain
$F_1\Delta_1(E)\otimes M_E\simeq F_2\Delta_2(E)\otimes N_E$ where
$M_E,N_E\in k\mproj$ and $\rank_k M_E=\rank_k N_E=\dim_K E$.
Since $\End_B(F_1\Delta_1(E))=k$, we deduce that
$\Hom_B(F_1\Delta_1(E),F_2\Delta_2(E))$ is an invertible $k$-module
and since
$$F_1\Delta_1(E)\otimes \Hom_B(F_1\Delta_1(E),F_2\Delta_2(E))\iso
F_2\Delta_2(E),$$
we obtain $F_2\Delta_2(E)\iso U_E\otimes F_1\Delta_1(E)$
for some $U_E\in\Pic(k)$.
The result follows now from Corollary \ref{unicityfc}.
\end{proof}

\begin{rem}
Let us give a variant of Theorem \ref{unicityfromorder}.
Let $\CC_1$ be a $1$-faithful highest weight cover of $B$ with associated
order $\le_1$ on $\Irr(KB)$.
Let $\le'$ be an order on $\Irr(KB)$ and $\{S'(E)\}_{E\in\Irr(KB)}$
be a set of $B$-modules such that
given $J'\subset I'$ coideals of $\Irr(KB)$ for $\le'$
such that no two distinct elements of $I'\setminus J'$ are comparable for
$\le'$, we have
$$\left((KB)^{I'}\cap B\right)/\left((KB)^{J'}\cap B\right)\simeq
\bigoplus_{E\in I'\setminus J'}S'(E)\otimes M_E$$
for some $M_E\in k\mproj$ with $\rank_k M_E=\dim_K E$.
Assume $\le_1$ is a refinement of $\le'$.
Then, given $E\in\Irr(KB)$, we have
$S'(E)\iso F_1\Delta_1(E)\otimes U_E$ for some $U_E\in\Pic(k)$.

\smallskip
In particular, if $\CC_2$ is a $1$-faithful highest weight cover of $B$ with
associated order $\le_2$ and if $\le_2$ is a refinement of $\le'$, then
$\CC_1$ and $\CC_2$ are equivalent highest weight covers.
\end{rem}

\begin{rem}
It would interesting to investigate when two $1$-faithful highest
covers are derived equivalent (cf Conjecture \ref{derequiv} for the case of
Cherednik algebras). This might be achieved through perverse
equivalences (cf \cite{ChRou} and \cite[\S 2.6]{Rou3}).
\end{rem}

\section{Cherednik's rational algebra}
\label{sectionCherednik}
We refer to \cite{Rou2} for a survey of the representation theory of
rational Cherednik algebras.

\subsection{Category $\CO$}
\subsubsection{}
\label{Cherednik}
Given $H\in\CA$, let $\alpha_H\in V^*$ with $H=\ker\alpha_H$ and let
$v_H\in V$ such that $\BC v_H$ is a $W_H$-stable complement to $H$.

The rational Cherednik algebra $\BA$ is the quotient of
$\BC[\{\Bh_u\}_{u\in U}]\otimes_\BC T(V\oplus V^*)\rtimes W$ by the relations
$$[\xi,\eta]=0 \text{ for } \xi,\eta\in V,\ \
[x,y]=0 \text{ for } x,y\in V^*$$
$$[\xi,x]=\langle \xi,x\rangle+\sum_{H\in\CA}
\frac{\langle\xi,\alpha_H\rangle\langle v_H,x\rangle}
     {\langle v_H,\alpha_H\rangle} \gamma_H$$
where
$$\gamma_H=\sum_{w\in W_H-\{1\}}
\left(\sum_{j=0}^{e_H-1}\det(w)^{-j}(\Bh_{H,j}-\Bh_{H,j-1})\right)w.$$

\begin{rem}
From the definition in \cite[\S 3.1]{GGOR} one gets to the notations here
by putting $h_{H,j}=-k_{H,-j}$ (here, we allow the possibility
$h_{H,0}\not=0$ to make twists by linear characters of $W$ more natural).

From the definitions in \cite[p.251]{EtGi}, $W$ being a finite Coxeter group,
one puts
$h_{H,0}=0$ and $h_{H,1}=c_{\alpha}$ for $H$ the kernel of the root
$\alpha$.
\end{rem}

\subsubsection{}
\label{secChhw}
Let $k'$ be a local commutative noetherian $\BC[\{\Bh_u\}]$-algebra with
residue field $k$.

Let $\CO$ be the category of finitely generated $k'\BA$-modules that are
locally nilpotent for $S(V)$.
Given $E\in\Irr(W)$, we put $\Delta(E)=k'\BA\otimes_{S(V)\rtimes W}E$
and we denote by $\nabla(E)$ the submodule of
$k'\Hom_{S(V^*)\rtimes W}(\BA,E)$ of elements that are locally finite for
$S(V)$.
Let $\Delta=\{\Delta(E)\}_{E\in\Irr(W)}$. We define an order on $\Irr(W)$
by $\chi>\chi'$ if $c_{\chi'}-c_{\chi}\in\BZ_{>0}$.

\begin{thm}
$(\CO,\Delta)$ is a highest weight category with 
costandard objects the $\nabla(E)$'s.
\end{thm}

\begin{proof}
We know that $\CO\simeq R\mMod$ for some finite projective $k'$-algebra
$R$ \cite[Corollary 2.8]{GGOR}. By Theorem \ref{carhwcat}, it suffices
to check the highest weight category property for $k\CO$: this is
given by \cite[Theorem 2.19]{GGOR}.
\end{proof}

\subsection{Covers of Hecke algebras}
\subsubsection{}
\label{KZ}

Let $\Gm$ be a maximal ideal of $\BC[\{\Bh_u\}]$ and $k'$ be the completion
at $\Gm$. We view $k'$ as a $\Bk$-algebra via
$\Bq_u\mapsto e^{2i\pi \Bh_u}$.

Let $\hat{\Gm}$ be the maximal ideal of $k'$ and
$k=k'/\hat{\Gm}$.
Let $h_\cdot=\{h_u\}\in k^U$ be the image of $\Bh$.
Let $\tGamma$ be the
subgroup of $k$ generated by $\BZ$ and the $h_u$'s.
We have an exact sequence
$$0\to \BZ\to\tGamma\xrightarrow{x\mapsto e^{2i\pi x}} e^{2i\pi\tGamma}\to 0$$
and we are in the setting of \S \ref{choice}, where we choose the coarsest
order.  In particular, the order on $\Irr(W)$ introduced in
\S \ref{secChhw} is the same as the one defined in \S \ref{order}.

There is a functor $\KZ:\CO\to k'\BH\mMod$ \cite[\S 5.3]{GGOR}
(note that in
the definition of the Hecke algebra in \cite[\S 5.2.5]{GGOR},
one should read $e^{-2i\pi\Bk_{H,j}}$ instead of $e^{2i\pi\Bk_{H,j}}$).
By \cite[\S 5.3, 5.4]{GGOR},
there is a projective object $P_{\KZ}$ of $\CO$ and an isomorphism
$k'\BH\iso\End_\CO(P_{\KZ})$ such that the functor $\KZ$ is isomorphic
to $\Hom_\CO(P_{\KZ},-)$.

\begin{thm}
\label{Cherednikcover}
$(k\CO,\KZ)$ is a highest weight cover of $k\BH$.

Assume
$x_{H,j}\not=x_{H,j'}$ for all $H\in\CA$ and $j\not=j'$. Then,
$(\CO,\KZ)$ is a $1$-faithful highest weight cover of $k'\BH$.
\end{thm}

\begin{proof}
The first statement is \cite[Theorem 5.16]{GGOR}. The second
statement follows, via Proposition \ref{equivDelta}, from
\cite[Proposition 5.9]{GGOR}.
\end{proof}

\begin{prop}
\label{ssCherednik}
Assume $\Gamma_{tor}=1$. Then, $k\CO$ and $k\BH$
are semi-simple.
\end{prop}

\begin{proof}
The semi-simplicities of $k\CO$ and of $k\BH$ are equivalent
(cf Theorem \ref{Cherednikcover}).
The algebra $k\BH$ depends only on the $h_u$'s up to shifts
by integers. So, in order to prove that $k\BH$ is semi-simple,
we can assume that the restriction of $t\mapsto e^{2i\pi t}$ to
the subgroup $\Gamma_0$ of $\BC$ generated by the $h_u$'s gives
an isomorphism $\Gamma_0\iso\Gamma$. Then, given
$\chi,\chi'\in\Irr(W)$, we have $c_\chi-c_{\chi'}\in\BZ$ if and only
if $c_\chi=c_{\chi'}$. In particular, no two distinct elements of
$\Irr(W)$ are comparable. So, $\CO$ is semi-simple and $k\BH$ as
well.
\end{proof}

\begin{proof}[Proof of Theorem \ref{ss}]
Without loss of generality, we may assume that $k$ has finite transcendence
degree over $\BQ$. Then, there is an embedding of $k$ in $\BC$ and
we can assume $k=\BC$. Now, the result follows from Proposition
\ref{ssCherednik}.
\end{proof}

\subsubsection{}
We denote by $\CO(h_\cdot)$ the category $k\CO$.

From Theorems \ref{unicityfromorder} and \ref{Cherednikcover}, 
we deduce a translation principle for category $\CO$:
\begin{thm}
\label{translation}
Assume $x_{H,j}\not=x_{H,j'}$ for all $H\in\CA$ and $j\not=j'$.
Let $\tau\in \Gt_\BZ$ and assume the order on $\Irr(W)$
defined by $h_\cdot$ is the
same as the one defined by $h_\cdot+\tau$.
Then, there is an equivalence $\CO(h_\cdot)\iso \CO(h_\cdot+\tau)$
of quasi-hereditary covers of $k\BH$.
\end{thm}

It would be interesting to describe precisely which $\tau$'s satisfy the
assumptions of the Theorem.

\begin{conj}
\label{derequiv}
Given any $\tau\in\Gt_\BZ$, then
$D^b(\CO(h_\cdot))\iso D^b(\CO(h_\cdot+\tau))$.
\end{conj}

\begin{rem}
Let $\kappa\in\BQ_{>0}$ with $\kappa\not\in(\frac{1}{e_H}\BZ)\setminus\BZ$
for all $H\in\CA$.
Assume $h_{H,j}=0$ for all $j\not=0$ and
$h_{H,0}=\kappa$, for all $H$. Let $\tau$ be given by
$\tau_{H,j}=0$ for $j\not=0$ and $\tau_{H,0}=1$.
Then $\tau$ satisfies the
assumption of the Theorem, \ie, the order defined by $h_\cdot$ is the same
as the one defined by $\tau+h_\cdot$.

\medskip
 We conjecture that, for general $W$,
the shift functor associated to $\zeta$ a linear character of $W$ gives an
equivalence if $h_\cdot$ and $h_\cdot+\tau$ define the same order on $\Irr(W)$,
where $\tau$ is the element corresponding to $\zeta$. Note that shift functors
are compatible with the KZ functor, hence when they are equivalences, they are
equivalences of highest weight covers of the Hecke algebra as in
Theorem \ref{translation}.

When $W$ has type $A_{n-1}$, Gordon and Stafford proved that the shift
functor is an equivalence (parameter ${\not\in}\ \frac{1}{2}+\BZ$)
\cite[Proposition 3.16]{GoSt1}.
\end{rem}

Note that equivalences arise also from twists \cite[\S 5.4.1]{GGOR}:

\begin{prop}
\label{equivtwist}
Let $\zeta\in W^\wedge$. 
We have an equivalence $\CO(h_\cdot)\iso \CO(\theta_\xi(h_\cdot))$ compatible, via
$\KZ$, with the isomorphism $\theta_\xi:\BH(\exp h_\cdot)\iso
 \BH(\exp\theta_\xi(h_\cdot))$.
\end{prop}

\subsubsection{}
We show now that Hecke algebras do not change, up to isomorphism of
$\BC$-algebras, by field automorphisms acting on parameters. As a consequence,
we show that category $\CO$ doesn't change if $h_\cdot$ is rescaled by a
positive integer, as long as the denominators do not change.

\smallskip
We fix $K_0$ be a subfield of $\BC$ such that the reflection representation
$V$ of $W$ is defined over $K_0$.

\begin{prop}
\label{HeckeoverK}
Let $q_\cdot\in \BT(\BC)$ with finite order.
Then, there exists a $K_0$-algebra $A$ and an isomorphism of $\BC$-algebras
$\BC\otimes_{K_0} A\simeq \BH(q_\cdot)$.
\end{prop}

\begin{proof}
Let $h_\cdot\in\BQ^U$ such that $q_\cdot=e^{2i\pi h_\cdot}$. Consider the
category $\CO_{K_0}$ for the rational Cherednik algebra defined over $K_0$, with
parameter $h_\cdot$. The simple objects of
$\CO_{K_0}$ remain simple in $\BC\otimes_{K_0}\CO_{K_0}$, hence there
is a projective object $P_{\KZ,K_0}$ of $\CO_{K_0}$ such that
$\BC\otimes_{K_0} P_{\KZ,K_0}\simeq P_{\KZ}$.
Then, $A=\End_{\CO_{K_0}}(P_{\KZ,K_0})$ satisfies the requirement
of the Proposition.
\end{proof}

For Hecke algebras, the next result, which is an immediate consequence
of Proposition \ref{HeckeoverK}, answers positively (over $\BC$)
a problem raised by Radha Kessar.
In type $A$, the result is due to Chuang and Miyachi \cite{ChMi}. Note that
their result covers also fields of positive characteristic.

\begin{thm}
\label{HeckeGalois}
Let $q_\cdot\in \BT(\BC)$ with finite order and let $\sigma$ be an automorphism
of $K_0(\{q_u\})/K_0$.

Then, we have an isomorphism of $\BC$-algebras:
$\BH(\sigma(q_\cdot))\simeq \BH(q_\cdot)$.
\end{thm}

\begin{rem}
\label{liftGalois}
The previous two results can be lifted. We use the notations of the proof
of Proposition \ref{HeckeoverK}.

Let $\Gm_0$ be the maximal ideal of $K_0[\{\Bh_u\}]$ generated by
the $\Bh_u-h_u$ and let $k_0$ be the completion at $\Gm_0$. Let
$k_{q_\cdot}$ be the completion of $\BC[\{\Bq_u^{\pm 1}\}]$ at
the maximal ideal generated by the $\Bq_u-q_u$.
Let $R=\BC[[\{X_u\}]]$ and consider the morphisms of algebras
$k_0\to R,\ \Bh_u-h_u\mapsto X_u$ and
$k_{q_\cdot}\to R,\ \Bq_u-q_u\mapsto e^{2i\pi X_u}$.
As in
Proposition \ref{HeckeoverK}, one shows there is a $K_0[\{\Bh_u\}]$-algebra
$A_0$ and an isomorphism of $R$-algebras $R\otimes_{k_0}A_0\simeq R
\otimes_{k_{q_\cdot}}k_{q_\cdot}\BH$.

\smallskip
Consider now the setting of Theorem \ref{HeckeGalois}.
We have an isomorphism of $R$-algebras
$R\otimes_{k_{\sigma(q_\cdot)}}k_{\sigma(q_\cdot)}\BH\simeq
 R \otimes_{k_{q_\cdot}}k_{q_\cdot}\BH$.
\end{rem}

\smallskip
\begin{thm}
\label{denominator}
Let $q_\cdot\in \BT(\BC)$ with finite order and let
$r\in\BZ_{>0}$ such that there is an automorphism
$\sigma\in K_0(\{q_u\})/K_0$ with $\sigma(q_\cdot)=q_\cdot^r$.
Assume $x_{H,j}\not=x_{H,j'}$ for all $H\in\CA$ and $j\not=j'$.

Then, there is an equivalence $\CO(h_\cdot)\iso \CO(rh_\cdot)$, which
identifies highest weight covers of $\BH(\sigma(q_\cdot))\simeq\BH(q_\cdot)$.
\end{thm}

\begin{proof}
The order on $\Irr(W)$ induced by $rh_\cdot$ is the same as the order induced
by $h_\cdot$. So, via the isomorphism of Remark \ref{liftGalois},
$\CO(h_\cdot)$  and $\CO(rh_\cdot)$ deform to $1$-faithful highest weight
covers of the same algebra (Theorem \ref{Cherednikcover}) and the result
follows from Theorem \ref{unicityfromorder}.
\end{proof}

Let us restate the previous Theorem in the case of Weyl groups and equal
parameters, where it takes a simpler form.

\begin{cor}
Assume $W$ is a Weyl group, $h_{u,1}=0$,
$h_{u,0}=h$ is constant and $h\in (\frac{1}{d}\BZ)
\setminus(\frac{1}{2}+\BZ)$ for some $d\in\BZ_{>0}$.
Given $r\in\BZ_{>0}$ prime to $d$, there is an equivalence
 $\CO(h)\iso \CO(rh)$, which
identifies highest weight covers of $\BH(q^r)\simeq\BH(q)$.
\end{cor}

Finally, let us relate characters. Define
$$\mathbf{eu}=\sum_{b\in\CB}b^\vee b+\sum_{H\in\CA}\sum_{j=1}^{e_H-1}
\sum_{w\in W_H}(h_{H,j}-h_{H,0})\det(w)^{-j}w$$
where $\CB$ is a basis of $V$ and $\{b^\vee\}_{b\in\CB}$ is the dual basis
of $V^*$. Given $M\in\CO$ and $a\in\BC$, we denote by $M_a$ the
generalized $a$-eigenspace of $\mathbf{eu}$ on $M$, a finite dimensional
vector space. The character of $M$ is an element of
$\BZ[[t]]\cdot\Irr(W)$ given by
$\chi_M(w,t)=\Tr_M(w\cdot t^{\mathbf{eu}})\in\BC[[t]]$ (here, $w\in W$).
Given $E\in\Irr(W)$, one has
$\chi_{\Delta(E)}(w,t)=\frac{\Tr_E(w)t^{c'_E}}{\det_{V^*}(1-wt)}$
(cf e.g. \cite[\S 2.1]{EtCh}).

\medskip
The following result was conjectured by Etingof. It follows immediately 
from Theorem \ref{denominator}.

\begin{prop}
With the assumptions of Theorem \ref{denominator}, we have
$$\chi_{L_{rh_\cdot}(E)}(w,t)=
\frac{\det_{V^*}(1-wt^r)}{\det_{V^*}(1-wt)}
\chi_{L_{h_\cdot}(E)}(w,t^r).$$
In particular, $L_{rh_\cdot}(E)$ is finite-dimensional if and only if
$L_{h_\cdot}(E)$ is finite-dimensional and when this is the case, we have
$\dim L_{rh_\cdot}(E)=r^{\dim V}\dim L_{h_\cdot}(E)$.
\end{prop}

\subsubsection{}
We discuss now blocks of ``defect $1$'' and show their structure depends only
on their number of simple objects.

\smallskip
Given $d$ a positive integer, recall that a Brauer tree algebra
associated to a line with $d$ vertices (and exceptional multiplicity $1$)
is a $\BC$-algebra Morita-equivalent to the principal block
of the Hecke algebra of the symmetric group $S_d$ at parameter
$(q_0,q_1)=(e^{2i\pi/d},-1)$
(cf \cite[\S 4.18]{Ben} for a general definition).
Consider now 
$$\widetilde{\mathrm{Br}}_d=\End_{\mathrm{Br}_d}(\mathrm{Br}_d\oplus \BC)$$
where $\BC$ is the trivial representation of $\mathrm{Br}_d$. This is
a quasi-hereditary algebra whose module category is ubiquitous in 
rational representation theory. It occurs as perverse sheaves on
$\BP^d$ for the partition $\BA^0\coprod \BA^1\coprod\cdots\coprod\BA^d$.

\medskip
We assume here that the algebra $\BH$ is endowed with a symmetrizing form
$t$: here, $t$ in a linear form $\BH\to\Bk$ with $t(ab)=t(ba)$ for all
$a,b\in\BH$ and the pairing $\BH\times\BH\to \Bk,\ (a,b)\mapsto t(ab)$,
is perfect. This is well-known to exist for $W$ a finite Coxeter group 
(take $t(T_w)=\delta_{1w}$) and
it is known to exist for the infinite series $G(r,p,n)$ \cite{MalMat}.

\smallskip
%
Let $\Gn\subset\hat{\Gm}$ be a prime ideal such that $R=k'/\Gn$
is a discrete valuation ring. Denote by $\pi$ a uniformizing parameter for
$R$. Denote by $K$ the field of fractions of $R$. Its residue field is $k$.

Let $\CA$ be a block of $R\CO$. We assume $K\CA$ is semi-simple.
We denote by $\Irr_\CA(W)$ the set of $E\in\Irr(W)$ such that
$\Delta(E)\in \CA$. We denote by $B$ the block of $R\BH$ corresponding,
via the $\KZ$-functor, to $\CA$ \cite[Corollary 5.18]{GGOR}.
Given $\chi\in\Irr_\CA(W)$, we denote by $s_\chi\in R$
the Schur element of $\chi_K$: the primitive idempotent of
$Z(K\BH)$ corresponding to $\chi_K$ is 
$s_\chi^{-1}\sum_a \chi_K(a)a^\vee$, where $a$ runs over a basis of $\BH$
over $\Bk$ and $(a^\vee)$ is the dual basis.

\medskip
The following Theorem gives the structure of blocks with defect one.
Theorem \ref{categorydefect1} was known for $W$ of type $A_n$ in
case the order of $h$ in $\BC/\BZ$ is $n+1$ \cite[Theorem 1.4]{BerEtGi}.
When $W$ is a Coxeter group, the statement about $\BH$ in Theorem
\ref{categorydefect1} goes back to Geck \cite[Theorem 9.6]{Ge1} (in the case of
equal parameters, but the proof applies to unequal parameters as well) and
we follow part of his proof.

\begin{thm}
\label{categorydefect1}
Let $d=|\Irr_\CA(W)|$.
Assume for every $\chi\in\Irr_\CA(W)$, we have
$\pi^{-1}s_\chi\in R^\times$ (``defect $1$'').
Then, $>$ is a total order on $\Irr_\CA(W)$, $kB$ is Morita equivalent to
$\mathrm{Br}_d$ and
$\CA$ is equivalent to $\widetilde{\mathrm{Br}}_d\mMod$.
In particular, if $\chi_1<\cdots<\chi_d$ are the elements of $\Irr_\CA(W)$,
then, for $n=1,\ldots,d$, we have
$$[L(\chi_n)]=\sum_{i=1}^n (-1)^{i+n}[\Delta(\chi_i)].$$
\end{thm}

\begin{proof}
Given $E\in\Irr_\CA(W)$, we denote by $L(E)$, $\Delta(E)$ and
$P(E)$ the corresponding simple, standard and projective objects of
$k\CA$.
Let $\Irr_\CA(W)^0$ be the set of $E\in\Irr_\CA(W)$ such that $\KZ(L(E))\not=0$.

\smallskip
Brauer's theory of blocks of finite groups of defect $1$ carries to 
$R\BH$ (cf \cite[Propositions 9.1-9.4]{Ge1} for the case of Weyl groups)
and shows that 
\begin{itemize}
\item[(i)]
$[\Delta(E):L(F)]\in\{0,1\}$ for $E\in\Irr_\CA(W)$ and $F\in\Irr_\CA(W)^0$.
\item[(ii)]
Given $E\in \Irr_\CA(W)^0$, there is a unique $F\in\Irr_\CA(W)$ distinct from $E$
such that $[P(E)]=[\Delta(E)]+[\Delta(F)]$.
\item[(iii)]
Given $E\in\Irr_\CA(W)$, then $\KZ(\Delta(E))$ is uniserial.
\end{itemize}

Let $E\in\Irr_\CA(W)$. Let $E_1\not= E_2\in\Irr_\CA(W)$ distinct from $E$ and 
such that $L(E_1)$ and $L(E_2)$ are composition factors of $\Delta(E)$. 

Since $[P(E_1)]=[\Delta(E_1)]+[\Delta(E)]$, the reciprocity formula
shows that $[\Delta(E_2):L(E_1)]=0$.
We have $[P(E_2)]=[\Delta(E_2)]+[\Delta(E)]$ by the reciprocity formula,
so we have an exact sequence
$$0\to \Delta(E)\to P(E_2)\to\Delta(E_2)\to 0.$$
Let $\Omega L(E_2)$ be the kernel of a projective cover
$P(E_2)\to L(E_2)$. Let $M$ be the kernel of a surjective map
$\Delta(E_2)\to L(E_2)$. We have an exact sequence
$$0\to \Delta(E)\to \Omega L(E_2)\to M\to 0.$$
Since $\Hom(M,L(E_1))=0$ and $\Hom(\Delta(E),L(E_1))=0$, it follows that
$\Hom(\Omega L(E_2),L(E_1))=0$, hence $\Ext^1(L(E_2),L(E_1))=0$.
Similarly, one shows that $\Ext^1(L(E_1),L(E_2))=0$.

Let $N$ be the kernel of a surjective map $\Delta(E)\to L(E)$. We have
shown that $N$ is semi-simple. Since $\KZ(\Delta(E))$ is uniserial, we deduce
that $\KZ(N)$ is simple or $0$.
So, we have proven 
\begin{itemize}
\item[(iv)]
Given $E\in\Irr_\CA(W)$, there is at most one $F\in\Irr_\CA(W)^0$ distinct from
$E$ and such that $[\Delta(E):L(F)]\not=0$.
\end{itemize}

The decomposition matrix of $B$ has at most two non-zero entries in
each row and in each column. It follows that $kB$ is a Brauer tree algebra
associated to a line (cf \cite[Theorem 9.6]{Ge1}). In particular,
the order $>$ on $\Irr_\CA(W)$ is a total order. Also, 
there is a unique $E'\in\Irr(W)$ such that $\KZ(L(E'))=0$.
We have $P(E')=\Delta(E')$ and $\KZ(\Delta(E'))$ is a simple module. Via
an appropriate
identification of $kB\mMod$ with $\mathrm{Br}_d\mMod$, it corresponds
to the trivial module $\BC$. Since $k\CA\simeq\End_{k\BH}\left(
\bigoplus_{E\in\Irr_\CA(W)}\KZ(P(E))\right)\mMod$, it follows that
$k\CA\simeq\widetilde{\mathrm{Br}}_d\mMod$.
\end{proof}

Let us give a concrete application of the previous result. Assume there is
$r\in\BZ_{>0}$ such that for all $u$, we have $h_u=\frac{a_u}{r}$
for some $a_u\in\BZ$.
The Schur element $s_\chi$ is the specialization at $\Bq_u=\Bq^{a_u}$ of
the generic Schur element $\Bs_\chi$ of $\chi$, where $\Bq=e^{2i\pi\Bh}$ and
$\pi=\Bh-\frac{1}{r}$. The assumption ``$\pi^{-1}s_\chi\in R^\times$'' will be
satisfied if and only if the $r$-th cyclotomic polynomial in $\Bq$ (over
$K_0$) divides $\Bs_\chi$ exactly once.

In case $a_u=1$ for all $u$ and $W$ is a finite Coxeter group, then
the principal block satisfies the assumption if and only if
$\Phi_r(q)$ divides the Poincar\'e polynomial of $W$ exactly once. Note that
in such a case the other blocks either satisfy the assumption or are 
simple.

\medskip
We list now for each finite exceptional irreducible Coxeter group $W$
all simple finite dimensional representations in a block $\CA$ of
defect $1$ and provide their character. We assume $a_u=1$ for all $u$.
We denote by $\phi_{m,b}$ an irreducible representation of $W$ of dimension
$m$ whose first occurrence in $S(V)$ is in degree $b$. When we use this
notation,
there is a unique irreducible representation of $W$ with that property. For
example, $\phi_{1,0}=\BC$ is the trivial representation and $\phi_{\dim V,1}=V$.
Computations have been performed in \GAP, using the \CHEVIE\ package
\cite{CHEVIE}. The blocks are described in \cite[Appendix F]{GePfe}.

\begin{itemize}
\item[$F_4$]
\begin{itemize}
\item[$\bullet$] $h=1/12$, $L(\BC)=\phi_{1,0}$.
\item[$\bullet$] $h=1/8$, $L(\BC)=\phi_{1,0}+t\phi_{4,1}+t^2\phi_{1,0}$.
\end{itemize}

\medskip
\item[$H_3$]
\begin{itemize}
\item[$\bullet$] $h=1/10$, $L(\BC)=\phi_{1,0}$.
\item[$\bullet$] $h=1/6$, $L(\BC)=\phi_{1,0}+t\phi_{3,1}+t^2\phi_{1,0}$.
\end{itemize}

\medskip
\item[$H_4$]
\begin{itemize}
\item[$\bullet$] $h=1/30$, $L(\BC)=\phi_{1,0}$.
\item[$\bullet$] $h=1/20$, $L(\BC)=\phi_{1,0}+t\phi_{4,1}+t^2\phi_{1,0}$.
\item[$\bullet$] $h=1/15$,
$L(\BC)=\phi_{1,0}+t\phi_{4,1}+t^2(\phi_{1,0}+\phi_{9,2})
+t^3\phi_{4,1}+t^4\phi_{1,0}$ and
$L(\phi_{4,7})=t^2\phi_{4,7}$.
\item[$\bullet$] $h=1/12$,
$L(\BC)=\phi_{1,0}+t\phi_{4,1}+t^2(\phi_{1,0}+\phi_{9,2})+
t^3(\phi_{4,1}+\phi_{16,3})+t^4(\phi_{1,0}+\phi_{9,2})+t^5\phi_{4,1}+
\phi_{1,0}$.
\item[$\bullet$] $h=1/10$,
$L(\phi_{4,1})=t^3(\phi_{4,1}+t\phi_{1,0}+t^2\phi_{4,1})$.
\end{itemize}

\medskip
\item[$E_6$]
\begin{itemize}
\item[$\bullet$] $h=1/12$, $L(\BC)=\phi_{1,0}$.
\item[$\bullet$] $h=1/9$, $L(\BC)=\phi_{1,0}+t\phi_{6,1}+t^2\phi_{1,0}$.
\end{itemize}

\medskip
\item[$E_7$]
\begin{itemize}
\item[$\bullet$] $h=1/18$, $L(\BC)=\phi_{1,0}$.
\item[$\bullet$] $h=1/14$, $L(\BC)=\phi_{1,0}+t\phi_{7,1}+t^2\phi_{1,0}$.
\item[$\bullet$] $h=1/10$, $L(\phi_{7,1})=t^{9/2}(\phi_{7,1}+t\phi_{1,0}+
t^2\phi_{7,1})$.
\end{itemize}

\medskip
\item[$E_8$]
\begin{itemize}
\item[$\bullet$] $h=1/30$,  $L(\BC)=\phi_{1,0}$.
\item[$\bullet$] $h=1/24$, $L(\BC)=\phi_{1,0}+t\phi_{8,1}+t^2\phi_{1,0}$.
\item[$\bullet$] $h=1/20$,
$L(\BC)=\phi_{1,0}+t\phi_{8,1}+t^2(\phi_{1,0}+\phi_{35,2})
+t^3\phi_{8,1}+t^4\phi_{1,0}$.
\item[$\bullet$] $h=1/18$,
$L(\phi_{8,1})=t^{5/3}(\phi_{8,1}+t(\phi_{1,0}+\phi_{28,8})+t^2\phi_{8,1})$.
\item[$\bullet$] $h=1/15$,
$L(\BC)=\phi_{1,0}+t\phi_{8,1}+t^2(\phi_{1,0}+\phi_{35,2})+
t^3(\phi_{8,1}+\phi_{112,3})+t^4(\phi_{1,0}+\phi_{35,2}+\phi_{210,4})+
t^4(\phi_{8,1}+\phi_{112,3})+t^5(\phi_{1,0}+\phi_{35,2})+t^6\phi_{8,1}
+t^7\phi_{1,0}$ and
$L(\phi_{8,1})=t^2(\phi_{8,1}+t(\phi_{1,0}+\phi_{28,8}+\phi_{35,2})+
t^2(2\phi_{8,1}+\phi_{160,7})+t^3(\phi_{1,0}+\phi_{28,8}+\phi_{35,2})+
t^4\phi_{8,1})$.
\item[$\bullet$] $h=1/12$,
$L(\phi_{28,8})=t^5(\phi_{28,8}+t(\phi_{8,1}+\phi_{56,19}) +t^2\phi_{28,8})$.
\end{itemize}
\end{itemize}

\begin{rem}
Consider a block $\CA$ of defect $1$.
Then, $\CA$ has at most one finite-dimensional
simple module. If $\CA$ has a finite-dimensional simple module,
it is $L(E)$ where $c'_E$ is minimal and we have
$|\Irr_\CA(W)|\ge 1+\dim V$,
since $L(E)$ has a projective resolution over $\BC[V]$ of length
$|\Irr_\CA(W)|$.
\end{rem}

\begin{rem}
It would be interesting to see if $|\Irr_\CA(W)|\le 1+\dim V$ for
any block $\CA$ satisfying the assumption of
Theorem \ref{categorydefect1}. Also, in case
of equal parameters with order $e$ in $\BC/\BZ$, is it true that
$|\Irr_\CA(W)|\le e$?
\end{rem}

\section{Case $W=B_n(d)$}
\label{sectionBnd}
\subsection{Combinatorics}
\subsubsection{}
Let $W$ be the complex reflection group of type $B_n(d)$ (\ie,
$G(d,1,n)$) for some integers $n,d\ge 1$.
This is the subgroup of $\GL_n(\BC)$ of monomial matrices
whose non-zero entries are $d$-th roots of unity.
The subgroup of permutations matrices is the symmetric group
$\GS_n$. It is generated by the transpositions
$s_1=(1,2),\ldots,s_{n-1}=(n-1,n)$.
Let $s_0$ be the diagonal matrix with diagonal
coefficients $(e^{2i\pi/d},1,\ldots,1)$. Then,
$W$ is generated by $s_0,s_1,\ldots,s_{n-1}$.
We identify its subgroup of diagonal matrices with the group
of functions $\{1,\ldots,n\}\to \mu_d$, where $\mu_d$ is the
group of $d$-th roots of unity of $\BC$.
Let $H_i$ be the reflecting hyperplane of $s_i$.

\smallskip
A partition of $n$ is a non-increasing sequence (finite or infinite)
$\alpha=(\alpha_1\ge\alpha_2\ge\ldots)$
of
non-negative integers with sum $n$ and we write $|\alpha|=n$.
We identify two partitions that differ only by zeroes.
We denote by ${^t\alpha}$ the transposed partition.
We denote by $\CP(n)$ the set of partitions of $n$.

A multipartition of $n$ is a $d$-tuple of partitions $\lambda=
(\lambda^{(1)},\ldots, \lambda^{(d)})$ with $\sum_i |\lambda^{(i)}|=n$.
We denote by $l_r$ the largest integer such that
$\lambda^{(r)}_{l_r}\not=0$.
We put
$$I_\lambda(r)=\{\sum_{i=1}^{r-1}|\lambda^{(i)}|+1,
\sum_{i=1}^{r-1}|\lambda^{(i)}|+2,\ldots,
\sum_{i=1}^r|\lambda^{(i)}|\}.$$
Given $i,j\ge 1$, we put
$$b_{i,j}^{(r)}=\begin{cases}
({^t \lambda}^{(r)})_j-i & \text{ if } ({^t \lambda}^{(r)})_j>i\\
0 & \text{ otherwise}
\end{cases}
\text{ and }
d_{i,j}^{(r)}=\begin{cases}
\lambda_i^{(r)}-j & \text{ if } \lambda_i^{(r)}>j \\
0 & \text{ otherwise.}
\end{cases}
$$
We put
$\GS_\lambda=\GS_{I_\lambda(1)}\times\cdots\GS_{I_\lambda(d)}$ and
$B_\lambda(d)=\mu_d^{\{1,\ldots,n\}}\rtimes \GS_\lambda$.
We denote by $\CP(d,n)$ the set of multipartitions of $n$.

\medskip
Given $\alpha\in\CP(n)$, we denote by $\chi_\alpha$ the
corresponding irreducible character of $\GS_n$.
Given $\lambda\in\CP(d,n)$, we denote by $\chi_\lambda$ the
corresponding irreducible character of $B_n(d)$. Let us
recall its construction.
We denote by $\phi^{(r)}$ the one-dimensional character of 
$(\mu_d)^{I_\lambda(r)}\rtimes\GS_{I_\lambda(r)}$ whose restriction
to $(\mu_d)^{I_\lambda(r)}$ is $\det^{r-1}$ and whose restriction to
$\GS_{I_\lambda(r)}$ is trivial.
Then,
$$\chi_\lambda=\Ind_{B_\lambda(d)}^{B_n(d)}
(\phi^{(1)}\chi_{\lambda^{(1)}}\otimes\cdots\otimes
\phi^{(d)}\chi_{\lambda^{(d)}}).$$

\begin{lemma}
\label{calculres}
Let $\lambda\in\CP(d,n)$. Then, given $0\le l\le d-1$, we have
$$\frac{1}{\chi_\lambda(1)}\langle(\chi_\lambda)_{|\langle s_0\rangle},
{\det}^l\rangle=
\frac{|\lambda^{(l+1)}|}{n}$$
and
$$\frac{1}{\chi_\lambda(1)}\langle(\chi_\lambda)_{|\langle s_1\rangle},
\det\rangle=\frac{1}{2}+\frac{1}{n(n-1)}
\sum_r \sum_{i,j} (b^{(r)}_{i,j}-d^{(r)}_{i,j}).$$
\end{lemma}

\begin{proof}
By Frobenius reciprocity and Mackey's formula, we have
$$\Res_{\langle s_0\rangle}\chi_\lambda=\sum_{i=1}^d
\frac{|\lambda^{(i)}| (n-1)!}{\prod_{r=1}^d |\lambda^{(r)}|!}
\left(\prod_{r=1}^d\chi_{\lambda^{(r)}}(1)\right)\cdot {\det}^{i-1}$$
hence
$$\frac{1}{\chi_\lambda(1)}\Res_{\langle s_0\rangle}\chi_\lambda=
\frac{1}{n}\sum_{i=1}^d |\lambda^{(i)}|\cdot {\det}^{i-1}.$$

We have
\begin{multline*}
\Res_{\langle s_1\rangle}\chi_\lambda=
\sum_{1\le r\le d,|\lambda^{(r)}|>1}
\frac{|\lambda^{(r)}|(|\lambda^{(r)}|-1)(n-2)!}{\prod_{i=1}^d |\lambda^{(i)}|!}
 \left(\prod_{1\le i\le d,i\not=r}\chi_{\lambda^{(i)}}(1)
\cdot \Res_{\GS_2}\chi_{\lambda^{(r)}}\right)+\\
+\sum_{r=1}^d \frac{|\lambda^{(r)}|(n-|\lambda^{(r)}|)(n-2)!}{2
\prod_{i=1}^d |\lambda^{(i)}|!}
\left(\prod_{i=1}^d\chi_{\lambda^{(i)}}(1)\right)(1+\det)
\end{multline*}
hence
\begin{multline*}
\frac{1}{\chi_\lambda(1)}\Res_{\langle s_1\rangle}\chi_\lambda=
\frac{1}{n(n-1)}
\Biggl(\sum_{1\le r\le d,|\lambda^{(r)}|>1}
\frac{|\lambda^{(r)}|(|\lambda^{(r)}|-1)}{\chi_{\lambda^{(r)}}(1)}
\cdot \Res_{\GS_2}\chi_{\lambda^{(r)}}+\\
+\sum_{r=1}^d\frac{|\lambda^{(r)}|(n-|\lambda^{(r)}|)}{2}\cdot (1+\det)
\Biggr).
\end{multline*}
Now, we have (cf Remark 3.3 and \cite[Theorem 10.5.2]{GePfe} for the generic degrees)
$$\frac{1}{\chi_{\lambda^{(r)}}(1)}
\langle \Res_{\GS_2}\chi_{\lambda^{(r)}},\det\rangle=
\frac{1}{2}+\frac{1}{|\lambda^{(r)}|(|\lambda^{(r)}|-1)}
\sum_{i,j}(b_{i,j}^{(r)}-d_{i,j}^{(r)})$$
and the second result follows.
\end{proof}

\subsubsection{}
Assume $d\not=1$ and $n\not=1$.
The braid group $B_W$ has generators $\sigma_0,\sigma_1,\ldots,
\sigma_{n-1}$ and relations \cite[Theorem 2.26]{BrMaRou}
$$\sigma_i\sigma_j=\sigma_j\sigma_i \textrm{ if } |i-j|>1,\
\sigma_0\sigma_1\sigma_0\sigma_1=\sigma_1\sigma_0\sigma_1\sigma_0
\textrm{ and } 
\sigma_i\sigma_{i+1}\sigma_i=\sigma_{i+1}\sigma_i\sigma_{i+1}
\textrm{ for } i\ge 1.$$
The canonical morphism $B_W\to W$ is given by
$\sigma_i\mapsto s_i$.

Put $\Bx_i=\Bx_{H_0,i}$, $\Bx=\Bx_{H_1,0}$ and put $\Bx_{H_1,1}=-1$.
Similarly, we will write
$h_i=h_{H_0,i}$, $h=h_{H_1,0}$ and assume $h_{H_1,1}=0$.

The Hecke algebra $\BH$ is the quotient of $\BZ[\Bq^{\pm 1},\Bx_0^{\pm 1},
\ldots,\Bx_{d-1}^{\pm 1}][B_W]$ by the ideal
generated by $(\sigma_0-\Bx_0)(\sigma_0-\Bx_1)\cdots(\sigma_0-\Bx_{d-1})$ and
$(\sigma_i-\Bq)(\sigma_i+1)$ for $1\le i\le n-1$ (this differs from the
algebra $\BH$ of \S\ref{secHecke} since we have already specialized
$\Bx_{H_1,1}$ to $-1$).

\bigskip
When $d=1$, then $B_W$ has generators $\sigma_1,\ldots,
\sigma_{n-1}$ and relations
$$\sigma_i\sigma_j=\sigma_j\sigma_i \textrm{ if } |i-j|>1
\textrm{ and } 
\sigma_i\sigma_{i+1}\sigma_i=\sigma_{i+1}\sigma_i\sigma_{i+1}\textrm{ for }
i\ge 1.$$
The canonical morphism $B_W\to W$ is given by $\sigma_i\mapsto s_i$.

Put $\Bx=\Bx_{H_1,0}$ and assume $\Bx_{H_1,1}=-1$. Similarly, let 
$h=h_{H_1,0}$ and assume $h_{H_1,1}=0$.

The Hecke algebra $\BH$ is the quotient of $\BZ[\Bq^{\pm 1}][B_W]$ by the ideal
generated by $(\sigma_i-\Bq)(\sigma_i+1)$ for $1\le i\le n-1$.

\bigskip
When $n=1$, then $B_W$ is an infinite cyclic group with one generator
$\sigma_0$. The canonical morphism $B_W\to W$ is given by
$\sigma_0\mapsto s_0$.

Put $\Bx_i=\Bx_{H_0,i}$ and $h_i=h_{H_0,i}$.

The Hecke algebra $\BH$ is the quotient of $\BZ[\Bx_0^{\pm 1},\ldots,
\Bx_{d-1}^{\pm 1}][B_W]$ by the ideal
generated by $(\sigma_0-\Bx_0)(\sigma_0-\Bx_1)\cdots(\sigma_0-\Bx_{d-1})$.

\bigskip
We denote by $T_i$ the image of $\sigma_i$ in $\BH$.
Note that $\BQ(\Bq,\Bx_0,\ldots,\Bx_{d-1})\BH$ is split semi-simple \cite{ArKo}.

\smallskip
From Lemma \ref{calculres}, we obtain
\begin{prop}
Let $\lambda\in\CP(d,n)$.
We have
$$c_{\chi_\lambda}=
d\sum_{2\le r\le d} |\lambda^{(r)}|(h_{r-1}-h_0)-
d\left(\frac{n(n-1)}{2}+\sum_{r,i,j}(b_{i,j}^{(r)}-d_{i,j}^{(r)})
\right)h.$$
\end{prop}

We put the dominance order $\unlhd$ on $\CP(d,n)$~: 
$\lambda\unlhd\mu$ if
$$\sum_{i=1}^{r-1} |\lambda^{(i)}|+ \sum_{j=1}^s |\lambda_j^{(r)}|\le
\sum_{i=1}^{r-1} |\mu^{(i)}|+ \sum_{j=1}^s |\mu_j^{(r)}|$$
for all $1\le r\le d$ and $s\ge 0$.

\begin{lemma}
\label{minimalorder}
Let $\lambda,\mu\in\CP(d,n)$. Then, $\lambda\lhd\mu$ and there is
no $\lambda'\in\CP(d,n)$ with $\lambda\lhd\lambda'\lhd\mu$ if and only
if one (or more) or the following holds~:
\begin{itemize}
\item[(a)] there is $s<d$ with 
\begin{itemize}
\item $\mu^{(r)}=\lambda^{(r)}$ for $r\not=s,s+1$
\item $\mu^{(s)}=(\lambda^{(s)}_1,\ldots,\lambda^{(s)}_{l_s},1)$
\item $\mu^{(s+1)}_1=\lambda^{(s+1)}_1-1$ and
$\mu^{(s+1)}_j=\lambda^{(s+1)}_j$ for $j>1$.
\end{itemize}
\item[(b)] there are $s$ and $i$ with
\begin{itemize}
\item $\mu^{(r)}=\lambda^{(r)}$ for $r\not=s$,
\item $\mu^{(s)}_j=\lambda^{(s)}_j$ for $j\not=i,i+1$,
$\mu^{(s)}_i=\lambda^{(s)}_i+1$ and
$\mu^{(s)}_{i+1}=\lambda^{(s)}_{i+1}-1$.
\end{itemize}
\item[(c)] there are $s$ and $i<i'$ with
\begin{itemize}
\item $\mu^{(r)}=\lambda^{(r)}$ for $r\not=s$,
\item $\mu^{(s)}_j=\lambda^{(s)}_j$ for $j\not=i,i'$ and
$\mu^{(s)}_i-1=\mu^{(s)}_{i'}+1=\lambda^{(s)}_i=\lambda^{(s)}_{i'}$.
\end{itemize}
\end{itemize}
\end{lemma}

\begin{proof}
Assume $\lambda\lhd\mu$ and there is
no $\lambda'\in\CP(d,n)$ with $\lambda\lhd\lambda'\lhd\mu$.
Take $s$ minimal such that $\lambda^{(s)}\not=\mu^{(s)}$.

Assume first that
$|\lambda^{(s)}|<|\mu^{(s)}|$. We denote by $m_s$ the largest integer such
that $\mu^{(s)}_{m_s}\not=0$.
If $\mu^{(s)}_{m_s}\not=1$, then
$\lambda\lhd \nu\lhd\mu$, where
$\nu^{(r)}=\mu^{(r)}$ for $r\not=s$ and
$\nu^{(s)}=(\mu^{(s)}_1,\ldots,\mu^{(s)}_{m_s-1},\mu^{(s)}_{m_s}-1,1)$, and
this is a contradiction. So, $\mu^{(s)}_{m_s}=1$.
Let $\xi\in\CP(d,n)$ be given by
$\xi^{(r)}=\mu^{(r)}$ for $r\not=s,s+1$,
$\xi^{(s)}=(\mu^{(s)}_1,\ldots,\mu^{(s)}_{m_s-1})$ and
$\xi^{(s+1)}=(\mu^{(s+1)}_1+1,\mu^{(s+1)}_2,\mu^{(s+1)}_3,\ldots)$.
Then,
$\lambda\unlhd\xi\lhd\mu$. So, $\lambda=\xi$ and we are in the case (a).

Assume now $|\lambda^{(s)}|=|\mu^{(s)}|$. Let
$\xi=(\mu^{(1)},\ldots,\mu^{(s)},\lambda^{(s+1)},\ldots,\lambda^{(d)})$.
Then, $\lambda\lhd\xi\unlhd\mu$, hence $\xi=\mu$. Now,
it is a classical fact about partitions that (b) or (c) holds
(cf e.g. \cite[Theorem 1.4.10]{JamKe}).

The other implication is clear.
\end{proof}

\begin{prop}
\label{comporders}
Assume $h\le 0$ and $h_s-h_{s-1}\ge (1-n)h$ for $1\le s\le d-1$.

Let $\lambda,\mu\in\CP(d,n)$.
If $\lambda\unlhd\mu$, then $c_{\chi_\lambda}\ge c_{\chi_\mu}$.
\end{prop}

\begin{proof}
It is enough to prove the Proposition in the case where $\lambda\not=\mu$ and
there is no $\lambda'$ with $\lambda\lhd\lambda'\lhd\mu$.
We use the description of Lemma \ref{minimalorder}.

Assume we are in case (a). Then,
$$c_{\chi_\lambda}-c_{\chi_\mu}=
d(h_{s}-h_{s-1})+dh(l_s+\lambda_1^{(s+1)}-1).$$
In case (b), we have
$$c_{\chi_\lambda}-c_{\chi_\mu}=
-dh(\mu_i^{(s)}-\mu_{i+1}^{(s)})$$
and in case (c), we have
$$c_{\chi_\lambda}-c_{\chi_\mu}=
-dh(i'-i+1).$$
The Proposition follows easily.
\end{proof}

\begin{rem}
One should compare the above order on $\CP(d,n)$ depending on $h$ and the
$h_i$'s to the order given by Jacon's $a$-function \cite[Definition 4.1]{Jac3}
and to the order defined by Yvonne \cite[\S 3.3]{Yv}.
\end{rem}

\subsection{The ``classical'' $q$-Schur algebras}
\label{classical}
\subsubsection{}
We recall here a generalization of Dipper and James' construction
(cf \cite{Do3}) of $q$-Schur algebras for type $A_{n-1}$ (case $d=1$ below).
As a first generalization, $q$-Schur algebras of type $B_n$
(case $d=2$ below) have been introduced
by Dipper, James, and Mathas \cite{DiJaMa1}, and Du and Scott \cite{DuSc2}.
The constructions have been then
extended by Dipper, James, and Mathas to the complex reflection
groups $B_n(d)$ \cite{DiJaMa2}.

\subsubsection{}

The subalgebra of $\BH$ generated by $T_1,\ldots,T_{n-1}$ is
the Hecke algebra of $\GS_n$, viewed as a Coxeter group with generating
set $(s_1=(1,2),\ldots,s_{n-1}=(n-1,n))$. Given $w=s_{i_1}\cdots
s_{i_r}\in\GS_n$, we put $T_w=T_{i_1}\cdots T_{i_r}$. We put
$L_i=\Bq^{1-i}T_{i-1}\cdots T_1T_0T_1\cdots T_{i-1}$.

Let $\lambda\in\CP(d,n)$. We put
$m_\lambda=\left(\sum_{w\in \GS_\lambda}T_w\right)
\left(\prod_{i=2}^d\prod_{j=1}^{a_i}(L_j-\Bx_i)\right)$,
where $a_i=|\lambda^{(1)}|+\cdots+|\lambda^{(i-1)}|$.

We put $M(\lambda)=m_\lambda\BH$, a right $\BH$-module, and
$P=\bigoplus_{\lambda\in\CP(d,n)}M(\lambda)$.
Let $\CS=\CS(d,n)=\End_{\BH^\opp}(P)^\opp$ (Dipper, James, and Mathas
consider a Morita equivalent algebra, where in the definition of $P$ the sum
is taken over all multicompositions of $n$).

\begin{thm}
\label{qSchurcover}
$(\CS,P)$ is a quasi-hereditary cover of $\BH$, for the order
given by the dominance order on $\CP(d,n)$.

Assume $k$ is a complete discrete valuation ring such that
$$(q+1)\prod_{i\not=j}(x_i-x_j)\in k^\times \textrm{ and }
\prod_{i=1}^n (1+q+\cdots+q^{i-1})\prod_{\substack{1\le i<j\le d\\
-n<r<n}}(q^rx_i-x_j)\not=0.$$

 Then
$(k\CS,kP)$ is a $1$-faithful quasi-hereditary cover of $k\BH$.
\end{thm}

\begin{proof}
The first assertion is known \cite[Theorems 4.14 and 5.3]{Mat2}.
The non-vanishing assumption is exactly the condition required to
ensure that $K\BH$ is split semi-simple \cite{Ar1}, where $K$ is the field
of fractions of $k$.
By \cite[Corollary 6.11 and Theorem 6.18]{Mat1}, given $T$ a tilting module
for $k\CS$, there is some
$k\BH$-module $M$ such that
$\Hom_{k\BH}(\Hom_{k\CS}(kP,k\CS),M)\simeq T$.
The second part follows now from
Propositions \ref{charff} and \ref{equivDelta}.
\end{proof}

\begin{rem}
In type $A$, these results are classical. Under the assumption that
$(1+q)(1+q+q^2)\not=0$ and $k$ is a field, then
$k\CS(1,n)$ is a $1$-faithful cover \cite[Theorem 3.8.1]{HeNa}.
See also \cite[\S 10]{Do4} for a different approach.
\end{rem}

We put $S(\lambda)=\Hom_{\CS}(P,\Delta(\lambda))$.

\subsection{Comparison}
\label{comparison}
In \S \ref{comparison}, we take $k,k'$ as in \S \ref{KZ}.

\subsubsection{}
The following result identifies category $\CO$ under certain assumptions.
\begin{thm}
\label{compBnd}
Assume
$(q+1)\prod_{i\not=j}(x_i-x_j)\not=0$.
Assume $h\le 0$ and $h_{s+1}-h_s\ge (1-n)h$ for $0\le s\le d-2$.

Then, $k\CO$ and
$k\CS\mMod$ are equivalent highest weight covers of $k\BH$~:
there is an equivalence
$k\CO\iso k\CS\mMod$ sending
the standard object associated to
$\chi\in\Irr(\GS_n)$ to the standard object associated to $\chi$.
\end{thm}

\begin{proof}
By Theorems
\ref{Cherednikcover} and \ref{qSchurcover},
$\CO$ and $k'\CS\mMod$ are $1$-faithful highest weight covers of $k'\CH$.
The order on irreducible
characters in $k'\CH$ coming from $k'\CS$ is a refinement of the one coming
from $\CO$, by Proposition \ref{comporders}. 
The Theorem follows now from Theorem \ref{unicityfromorder}.
\end{proof}

Note that, under the assumptions of the Theorem, $\CO$ and
$k'\CS\mMod$ are equivalent highest weight covers of $k'\BH$ as well.

\begin{rem}
Using Proposition \ref{equivtwist}, we obtain other parameter values for which
$k\CO$
is equivalent to $k\CS\mMod$ (for example, replacing $h$ by $-h$ in the
Theorem). The Theorem should hold without the assumption
$(q+1)\prod_{i\not=j}(x_i-x_j)\not=0$, but the methods developed here cannot
handle this general case.
\end{rem}

\begin{rem}
This suggests to look for a construction similar to that of \S \ref{classical} of
$q$-Schur algebras of type $B_n(d)$ for orders on $\CP(d,n)$ coming from
other choices of $h$ and $h_i$'s. Recent work of Gordon \cite{Go}
provides an order
based on the geometry of Hilbert schemes that is probably more relevant that
the orders used here.

It might be possible to produce
explicit ``perverse complexes'' and obtain the other $q$-Schur algebras
by perverse tilts (cf Conjecture \ref{derequiv}).
\end{rem}

\subsubsection{}
Let us restate the previous Theorem in the case $W=\GS_n$. In that case,
$\CS(1,n)$ is the $q$-Schur algebra of $\GS_n$, Morita equivalent
to a quotient of the quantum group $U_q(\gl_n)$. The following result
solves a conjecture of \cite[Remark 5.17]{GGOR} (under the assumption
$h\not\in\frac{1}{2}+\BZ$).

\begin{thm}
\label{mainsymm}
Assume $h\not\in\frac{1}{2}+\BZ$. Then, there is
an equivalence of highest weight categories
$k\CO\iso k\CS(1,n)\mMod$ sending the standard object associated to
$\chi\in\Irr(\GS_n)$ to the standard object associated to
$\begin{cases}
\chi & \textrm{ if } h\le 0 \\
\chi\otimes\det & \textrm{ if } h>0.
\end{cases}$
\end{thm}

This shows the characters of simple objects of $\CO$ are given by
canonical basis elements in the Fock space for $\hat{\Gsl}_r$, where
$r$ is the order of $k$ in $\BC/\BZ$, according to 
Varagnolo-Vasserot's proof \cite{VarVas1} of Leclerc-Thibon's conjecture
\cite{LeTh} (a generalization of Ariki's result \cite{Ar2} proving
Lascoux-Leclerc-Thibon's conjecture \cite{LaLeTh}).
Cf \S \ref{Uglov} for a conjectural generalization to the case $d>1$.

Gordon and Stafford \cite[Proposition 6.11]{GoSt2} deduce from this result
a description
of the maximal dimensional components of the characteristic cycle of the
simple objects (a cycle in $\Hilb^n\BC^2$). If these
characteristic cycles were equidimensional, they would thus be known and one
could deduce what are the support varieties in $\BC^{2n}/\GS_n$
of the simple objects in $\CO$.

\begin{rem}
One can expect to obtain a different proof of Theorem \ref{mainsymm} via
the work of Suzuki \cite{Su}, which relates representations of rational
Cherednik algebras of type $A$ with representations at negative level
of affine Lie algebras of type $A$.
\end{rem}

\smallskip
Note that an analog of Theorem \ref{mainsymm} has been proven by
Varagnolo and Vasserot for trigonometric (or elliptic) Cherednik algebras
\cite{VarVas2}.

\subsection{Orbit decomposition}
Let $s\in\{0,\ldots,d-1\}$ such that
$(q^ix_r-q^{i'}x_{r'})\in k^\times$ for
$0\le r<s$, $s\le r'<d$ and $0\le i,i'\le n$.
There is a bijection
\begin{align*}
\coprod_{m=0}^n\CP(s,m)\times\CP(d-s,n-m)&\xrightarrow[\sim]{\cup} \CP(d,n)\\
(\alpha^{(1)},\ldots,\alpha^{(s)}),(\beta^{(1)},\ldots,\beta^{(d-s)})
&\mapsto 
(\alpha^{(1)},\ldots,\alpha^{(s)},\beta^{(1)},\ldots,\beta^{(d-s)}).
\end{align*}

We write $\BH_{x_0,\ldots,x_{d-1}}(n)$ for the algebra $k\BH$ (which
depends further on $q$).

In \cite[Theorem 1.6]{DiMa}, Dipper and Mathas construct an
equivalence
$$F:\left(\bigoplus_{m=0}^n
\BH_{x_0,\ldots,x_{s-1}}(m)\otimes
\BH_{x_s,\ldots,x_{d-1}}(n-m)\right)\mMod\iso \BH_{x_0,\ldots,x_{d-1}}(n)\mMod$$
with the property that $F(S(\alpha)\otimes S(\beta))=S(\alpha\cup\beta)$
\cite[Proposition 4.11]{DiMa}.

\medskip
Assume we are in the setting of \S \ref{KZ}.
We write $\CO_{h_0,\ldots,h_{d-1}}(n)$ for the category $k\CO$.

\begin{thm}
Assume $(q+1)\prod_{i\not=j}(x_i-x_j)\not=0$.
Let $s\in\{0,\ldots,d-1\}$ such that
$q^ix_r\not=q^{i'}x_{r'}$ for
$0\le r<s$, $s\le r'<d$ and $0\le i,i'\le n$.
Then, there is an equivalence of highest weight categories
$$\left(\bigoplus_{m=0}^n
\CO_{h_0,\ldots,h_{s-1}}(m)\otimes
\CO_{h_s,\ldots,h_{d-1}}(n-m)\right)\mMod\iso
\CO_{h_0,\ldots,h_{d-1}}(n)\mMod.$$
It sends $\Delta(\alpha)\otimes\Delta(\beta)$ to $\Delta(\alpha\cup\beta)$
and it is compatible with $F$.
\end{thm}

\begin{proof}
Fix $m$ and consider $\alpha\in\CP(s,m)$ and $\beta\in\CP(d-s,n-m)$.
We have
\begin{multline*}
\frac{c_{\chi_\alpha}}{s}+\frac{c_{\chi_\beta}}{d-s}=
\sum_{2\le r\le s} |\alpha^{(r)}|(h_{r-1}-h_0)+
\sum_{s+1\le r\le d} |\beta^{(r-s)}|(h_{r-1}-h_s)-\\
-\left(\frac{m(m-1)+(n-m)(n-m-1)}{2}+\sum_{r,i,j}(b_{i,j}^{(r)}-d_{i,j}^{(r)})
\right)h
\end{multline*}
so
$$d(\frac{c_{\chi_\alpha}}{s}+\frac{c_{\chi_\beta}}{d-s})=
c_{\chi_{\alpha\cup\beta}}+d(n-m)(h_0-h_s)+m(n-m).$$
We deduce that if $\chi_\alpha\le\chi_{\alpha'}$ and
$\chi_\beta\le\chi_{\beta'}$, then $\chi_{\alpha\cup\beta}\le
\chi_{\alpha'\cup\beta'}$.

The result follows now from Theorems \ref{Cherednikcover} and
\ref{unicityfromorder}.
\end{proof}

\begin{rem}
The Theorem should hold without the assumption 
$(q+1)\prod_{i\not=j}(x_i-x_j)\not=0$.
\end{rem}

\begin{rem}
Note that this Theorem applies to more general $1$-faithful highest weight
covers (in particular, to the classical one, where we recover
\cite[Theorem 1.5]{DiMa}, with the additional assumption that
$(q+1)\prod_{i\not=j}(x_i-x_j)\in k^\times$).
\end{rem}

\begin{rem}
We put an equivalence relation on $\{0,1,\ldots,d-1\}$:~$r$ and $r'$ are
equivalent if there is $a\in\{-n,\ldots,n\}$ such that
$x_{r'}=q^ax_r$. Then, $\CO$ is equivalent to 
$$\bigoplus_{\substack{m:(\{0,\ldots,d-1\}/\sim)\to \BZ_{\ge 0}\\
\sum_I m(I)=n}}
\bigotimes_{I\in \{0,\ldots,d-1\}/\sim} \CO_{\{h_i\}_{i\in I}}(m(I)).$$
\end{rem}

\subsection{Uglov's higher level Fock spaces}
\label{Uglov}
Let $e>1$ be an integer and let $s_.=(s_0,\ldots,s_{d-1})\in\BZ^d$.
Let $h=\frac{1}{e}$ and $h_j=\frac{s_j}{e}-\frac{j}{d}$.
Uglov \cite{Ug} has introduced a $q$-deformed Fock space of level $d$
associated to the multicharge $s_.$, together with a standard and a canonical
basis, both parametrized by $d$-multipartitions.

Yvonne \cite{Yv} conjectured that, for suitable values of the $s_i$'s,
the multiplicities of simple modules in
standard modules for classical $q$-Schur algebras are equal to the
corresponding coefficients of the
transition matrix between the standard and the canonical basis.
He showed that this is compatible with the Jantzen sum formula.

Now, Theorem \ref{compBnd} shows that Yvonne's conjecture can be restated
for category $\CO$ and particular values of the $s_i$'s. We
conjecture that, for arbitrary $s_i$'s, the multiplicities of simple modules
in standard modules in $\CO$ are equal to the
corresponding coefficients of the
transition matrix between the standard and the canonical basis.
We also expect that the $q$-coefficients measure the level in the 
filtration induced by the Shapovalov form.
It should be possible to prove a sum formula for Cherednik algebras
and obtain a result similar to Yvonne's.

\begin{rem}
In order to prove the conjecture (in the case ``$(q+1)\prod_{i\not=j}(x_i-x_j)
\not=0$''), it would suffice to construct a (deformation of a) highest
weight cover of the Hecke algebra of a geometrical nature, where character
formulas can be computed.
\end{rem}


\begin{thebibliography}{[BrMaRou]}
\bibitem[Ar1]{Ar1} S.~Ariki,
	{\em On the semi-simplicity of the Hecke algebra of
	$(\BZ/rBZ)\wr \GS_n$}, J. Algebra {\bf 169} (1994), 216--225.
\bibitem[Ar2]{Ar2} S.~Ariki,
	{\em On the decomposition numbers of the Hecke algebra of $G(m,1,n)$},
	J. Math. Kyoto Univ. {\bf 36} (1996), 789--808.
\bibitem[ArKo]{ArKo} S.~Ariki and K.~Koike,
	{\em  A Hecke algebra of $(Z/rZ)\wr\GS_n$ and construction
	of its irreducible representations},
	Adv. Math. {\bf 106} (1994), 216--243.
\bibitem[Ben]{Ben} D.~Benson,
	``Representations and cohomology I'', Cambridge University Press,
	1991.
\bibitem[BerEtGi]{BerEtGi} Y.~Berest, P.~Etingof and V.~Ginzburg,
	{\em Finite-dimensional representations of rational Cherednik
	algebras}, Int. Math. Res. Not.  {\bf 19} (2003), 1053--1088.
\bibitem[BrMaMi]{BrMaMi} M.~Brou\'e, G.~Malle and J.~Michel,
	{\em Towards spetses. I},
	Transform. Groups {\bf 4} (1999), 157--218.
\bibitem[BrMaRou]{BrMaRou} M.~Brou\'e, G.~Malle and R.~Rouquier,
        {\em Complex reflection groups, braid groups, Hecke algebras},
        J. reine angew. Math. {\bf 500} (1998), 127--190.
\bibitem[BrMi]{BrMi} M.~Brou\'e and J.~Michel,
        {\em  Sur certains \'el\'ements r\'eguliers des groupes de Weyl
        et les vari\'et\'es de Deligne-Lusztig associ\'ees},
        in ``Finite reductive groups'', Birkh\"auser, 73--139, 1997.
\bibitem[CHEVIE]{CHEVIE} M.~Geck, G.~Hi{\ss}, F.~L\"ubeck, G.~Malle and
	G.~Pfeiffer,
	{\em {\sf CHEVIE} --- A system for computing and processing generic
	character tables},
	Appl. Algebra Engr. Comm. Comput. {\bf 7} (1996), 175--210.
\bibitem[ChMi]{ChMi} J.~Chuang and H.~Miyachi,
	{\em Runner removal Morita equivalences}, preprint, 2007.
\bibitem[ChRou]{ChRou} J.~Chuang and R.~Rouquier,
	{\em Calabi-Yau algebras and perverse Morita equivalences},
	in preparation.
\bibitem[CPS1]{CPS1} E.~Cline, B.~Parshall and L.~Scott,
	{\em  Finite-dimensional algebras and highest weight categories},
	J. Reine Angew. Math. {\bf 391} (1988), 85--99.
\bibitem[CPS2]{CPS2} E.~Cline, B.~Parshall and L.~Scott,
	{\em Integral and graded quasi-hereditary algebras, I},
	J. of Alg. {\bf 131} (1990), 126--160.
\bibitem[DiJaMa1]{DiJaMa1} R.~Dipper, G.~James and A.~Mathas,
	{\em The $(Q,q)$-Schur algebra},
	Proc. London Math. Soc. {\bf 77} (1998), 327--361.
\bibitem[DiJaMa2]{DiJaMa2} R.~Dipper, G.~James and A.~Mathas,
	{\em Cyclotomic $q$-Schur algebras},
	Math. Z. {\bf 229} (1998), 385--416.
\bibitem[DiMa]{DiMa} R.~Dipper and A.~Mathas,
	{\em Morita equivalences of Ariki-Koike algebras},
	 Math. Z. {\bf 240} (2002), 579--610.
\bibitem[Do1]{Do1} S.~Donkin,
	{\em On tilting modules for algebraic groups},
	Math. Z. {\bf 212} (1993), 39--60.
\bibitem[Do2]{Do2} S.~Donkin,
	{\em  On Schur algebras and related algebras. III.
	Integral representations},
	Math. Proc. Cambridge Philos. Soc. {\bf 116} (1994), 37--55.
\bibitem[Do3]{Do3} S.~Donkin,
	``The $q$-Schur algebra'',
	Cambridge University Press, 1998.
\bibitem[Do4]{Do4} S.~Donkin,
	{\em Tilting modules for algebraic groups and finite dimensional
	algebras}, in ``Handbook of tilting theory'', 215--257,
	Cambridge Univ. Press, 2007.
\bibitem[DuPaSc1]{DuPaSc1} J.~Du, B.~Parshall and L.~Scott,
	{\em Stratifying endomorphism algebras associated to Hecke algebras},
	J. Algebra {\bf 203} (1998), 169--210.
\bibitem[DuPaSc2]{DuPaSc2} J.~Du, B.~Parshall and L.~Scott,
	{\em Cells and $q$-Schur algebras},
	Transf. Groups {\bf 3} (1998), 33--49.
\bibitem[DuPaSc3]{DuPaSc3} J.~Du, B.~Parshall and L.~Scott,
	{\em  Quantum Weyl reciprocity and tilting modules},
	Comm. Math. Phys. {\bf 195} (1998), 321--352.
\bibitem[DuSc1]{DuSc1} J.~Du and L.~Scott,
	{\em Lusztig conjectures, old and new, I},
	J. Reine Angew. Math. {\bf 455} (1994), 141--182.
\bibitem[DuSc2]{DuSc2} J.~Du and L.~Scott,
	{\em The $q$-${\rm Schur}\sp 2$ algebra},
	Trans. Amer. Math. Soc. {\bf 352} (2000), 4325--4353.
\bibitem[DunOp]{DunOp} C.~Dunkl and E.~Opdam,
        {\em Dunkl operators for complex reflection groups},
	Proc. London Math. Soc. {\bf 86} (2003), 70--108.
\bibitem[EtCh]{EtCh}  P.~Etingof and T.~Chmutova,
	{\em On some representations of the rational Cherednik algebra},
	Representation Theory {\bf 7} (2003), 641--650.
\bibitem[EtGi]{EtGi}  P.~Etingof and V.~Ginzburg,
        {\em Symplectic reflection algebras, Calogero-Moser space and
        deformed Harish-Chandra homomorphism}, Inv. Math. {\bf 147}
        (2002), 243--348.
\bibitem[EtRa]{EtRa}  P.~Etingof and E.~Rains,
	{\em Central extensions of preprojective algebras, the quantum
	Heisenberg algebra, and $2$-dimensional complex reflection groups},
	J. Algebra {\bf 299} (2006), 570--588.
\bibitem[Ge1]{Ge1} M.~Geck,
	{\em Brauer trees of Hecke algebras},
	Comm. Algebra {\bf 20} (1992), 2937--2973.
\bibitem[Ge2]{Ge2} M.~Geck,
	{\em Kazhdan-Lusztig cells and decomposition numbers},
	Represent. Theory {\bf 2} (1998), 264--277.
\bibitem[Ge3]{Ge3} M.~Geck,
	{\em Kazhdan-Lusztig cells, $q$-Schur algebras and James' conjecture},
	J. London Math. Soc. {\bf 63} (2001), 336--352.
\bibitem[GePfe]{GePfe} M.~Geck and G.~Pfeiffer,
	``Characters of finite Coxeter groups and Iwahori-Hecke algebras'',
	Oxford University Press, 2000.
\bibitem[GeRou]{GeRou} M.~Geck and R.~Rouquier,
	{\em  Filtrations on projective modules for Iwahori-Hecke algebras},
	in ``Modular representation theory of finite groups'', de Gruyter,
	211--221, 2001.
\bibitem[GGOR]{GGOR}V.~Ginzburg, N.~Guay, E.~Opdam and R.~Rouquier,
        {\em On the category $\mathcal{O}$ for rational Cherednik algebras},
	Inventiones Math. {\bf 154} (2003), 617--651
\bibitem[Go]{Go} I.~Gordon,
	{\em  Quiver varieties, category $\mathcal{O}$ for rational Cherednik
	algebras, and Hecke algebras}, preprint math/0703150.
\bibitem[GoSt1]{GoSt1} I.~Gordon and J.T.~Stafford,
	{\em Rational Cherednik algebras and Hilbert schemes},
	Adv. Math. {\bf 198} (2005), 222--274.
\bibitem[GoSt2]{GoSt2} I.~Gordon and J.T.~Stafford,
	{\em Rational Cherednik algebras and Hilbert schemes, II:
	representations and sheaves},
	Duke Math. J. {\bf 132} (2006), 73--135.
\bibitem[HeNa]{HeNa} D.J.~Hemmer and D.K.~Nakano,
	{\em Specht filtrations for Hecke algebras of type A},
	J. London Math. Soc. {\bf 69}  (2004), 623--638.
\bibitem[Jac1]{Jac1} N.~Jacon,
	{\em  Sur les repr\'esentations modulaires des alg\`ebres de Hecke de
	type $D\sb n$},
	J. Algebra {\bf 274} (2004), 607--628.
\bibitem[Jac2]{Jac2} N.~Jacon,
	{\em On the parametrization of the simple modules for Ariki-Koike
	algebras at roots of unity},
	J. Math. Kyoto Univ. {\bf 44} (2004), 729--767.
\bibitem[Jac3]{Jac3} N.~Jacon,
	{\em Crystal graphs of higher level q-deformed Fock spaces,
	Lusztig a-values and Ariki-Koike algebras},
	preprint math.RT/0504267, to appear in Algebras and Representation
	Theory.
\bibitem[JamKe]{JamKe} G.~James and A.~Kerber,
	``The representation theory of the symmetric group'',
	Addison-Wesley, 1981.
\bibitem[K\"o]{Ko} S.~K\"onig,
	{\em  Quasi-hereditary orders},
	Manuscripta Math. {\bf 68} (1990), 417--433.
\bibitem[LaLeTh]{LaLeTh} A.~Lascoux, B.~Leclerc and J.-Y.~Thibon,
	{\em Hecke algebras at roots of unity and crystal bases of
	quantum affine algebras}, Comm. Math. Phys.
	{\bf 181} (1996), 205--263.
\bibitem[LeTh]{LeTh} B.~Leclerc and J.-Y.~Thibon,
	{\em  Canonical bases of $q$-deformed Fock spaces},
	Internat. Math. Res. Notices {\bf 9} (1996), 447--456.
\bibitem[Lu1]{Lu1} G.~Lusztig,
	{\em  Left cells in Weyl groups},
	in ``Lie group representations, I'', 99--111,
	Lecture Notes in Math. 1024, Springer, Berlin, 1983
\bibitem[Lu2]{Lu2} G.~Lusztig,
	``Hecke algebras with unequal parameters'',
	CRM Monograph Series, 18, Amer. Math. Soc., 2003.
\bibitem[Mal]{Mal} G.~Malle,
	{\em  On the rationality and fake degrees of characters of
	cyclotomic algebras},
	J. Math. Sci. Univ. Tokyo  {\bf 6}  (1999), 647--677.
\bibitem[MalMat]{MalMat} G.~Malle and A.~Mathas,
	{\em Symmetric cyclotomic Hecke algebras},
	J. Algebra {\bf 205} (1998), 275--293. 
\bibitem[Mat1]{Mat1} A.~Mathas,
	{\em Tilting modules for cyclotomic Schur algebras},
	J. Reine Angew. Math. {\bf 562} (2003), 137--169.
\bibitem[Mat2]{Mat2} A.~Mathas,
	{\em The representation theory of the Ariki-Koike and cyclotomic
	$q$-Schur algebras},
	in ``Representation theory of algebraic groups and quantum groups'',
	pp. 261--320, Adv. Stud. Pure Math., 40, Math. Soc. Japan, 2004.
\bibitem[Rou1]{Rou1} R.~Rouquier,
	{\em Familles et blocs d'alg\`ebres de Hecke},
	C. R. Acad. Sci. Paris {\bf 329} (1999), 1037--1042.
\bibitem[Rou2]{Rou2} R.~Rouquier,
	{\em Representations of rational Cherednik algebras},
	in ``Infinite-dimensional aspects of representation theory and
	applications'', 103--131, Amer. Math. Soc., 2005.
\bibitem[Rou3]{Rou3} R.~Rouquier,
	{\em Derived equivalences and finite dimensional algebras},
	in Proceedings ICM Madrid 2006, vol. II, 191--221, European
	Mathematical Society, 2006.
\bibitem[Su]{Su} T.~Suzuki,
	{\em Double affine Hecke algebras, conformal coinvariants and Kostka
	polynomials},
	 C. R. Math. Acad. Sci. Paris {\bf 343} (2006), 383--386.
\bibitem[Ug]{Ug} D.~Uglov,
	{\em Canonical bases of higher-level $q$-deformed Fock spaces
	and Kazhdan-Lusztig polynomials},
	in ``Physical combinatorics (Kyoto, 1999)'', pp. 249--299,
	Birkh\"auser, 2000. 
\bibitem[VarVas1]{VarVas1} M.~Varagnolo and E.~Vasserot,
	{\em  On the decomposition matrices of the quantized Schur algebra},
	Duke Math. J. {\bf 100} (1999), 267--297.
\bibitem[VarVas2]{VarVas2} M.~Varagnolo and E.~Vasserot,
	{\em  From double affine Hecke algebras to quantized affine Schur
	algebras}, Int. Math. Res. Not. (2004), {\bf 26}, 1299--1333.
\bibitem[Yv]{Yv} X.~Yvonne,
	{\em A conjecture for q-decomposition matrices of cyclotomic
	v-Schur algebras},
	J. Algebra {\bf 304} (2006), 419--456. 
\end{thebibliography}
\end{document}